\documentclass[aps,pre,onecolumn,11pt,notitlepage,nofootinbib]{amsart}
\usepackage[margin=0.8in]{geometry}
 \pdfoutput=1
\usepackage{setspace}
\usepackage{bm,multirow}
\usepackage[dvipsnames]{xcolor}
\usepackage{epsfig,graphics,amssymb,amsmath,amsthm,color,setspace}
\usepackage[normalem]{ulem}
\usepackage{graphicx,txfonts,listings}
\usepackage{booktabs}
\usepackage{array}
\usepackage[sfdefault=cmbr]{isomath}
\usepackage[colorinlistoftodos]{todonotes}

\newcommand{\ra}[1]{\renewcommand{\arraystretch}{#1}}

\author{William H. Mitchell}
\author{Henry G. Bell}
\author{Yoichiro Mori}
\author{Laurel Ohm}
\author{Daniel Spirn}

\begin{document}
\title{{A single-layer based numerical method for the slender body boundary value problem.}}
\date{\today}

\begin{abstract}
Fluid flows containing dilute or dense suspensions of thin fibers are widespread in biological and industrial processes. To describe the motion of a thin immersed fiber, or to describe the forces acting on it, it is convenient to work with one-dimensional fiber centerlines and force densities rather than two-dimensional surfaces and surface tractions. Slender body theories offer ways to model and simulate the motion of immersed fibers using only one-dimensional data. However, standard formulations can break down when the fiber surface comes close to intersecting itself or other fibers. In this paper we introduce a numerical method for a recently derived three-dimensional \emph{slender body boundary value problem} that can be stated entirely in terms of a one-dimensional distribution of forces on the centerline. The method is based on a new \emph{completed single-layer potential} formulation of fluid velocity which removes the nullspace associated with the unmodified single layer potential. We discretize the model and present numerical results demonstrating the good conditioning and improved performance of the method in the presence of near-intersections. To avoid the modeling and numerical choices involved with free ends, we consider closed fibers. 
\end{abstract}

\maketitle

\section{Introduction}
The use of small parameters to simplify difficult modeling and simulation problems is one of the outstanding successes of classical applied mathematics. In fluid mechanics, one important example of a small parameter is the aspect ratio of an immersed structure. Biological examples of slender immersed objects include microtubules inside cells  \cite{shelley2016dynamics,nazockdast2017fast} and cilia and flagella external to cells \cite{rodenborn2013propulsion,lauga2009hydrodynamics,cortez2005method,lighthill1976flagellar}; many industrial processes such as papermaking also rely on the properties of sparse or dense fiber suspensions \cite{hamalainen2011papermaking,petrie1999rheology}. It is very appealing to work with the one-dimensional centerlines of these thin structures rather than with their two-dimensional surfaces or three-dimensional volumes. In addition to the computational efficiency that comes with lowering the dimension of the problem, there are also important theoretical advantages.  For example, it is much simpler to formulate a model for the centerline density of forces on a fiber than it is to model fully two-dimensional surface tractions. 

An attempt to make physical sense of forces and velocities defined on the one-dimensional centerline instead of the two-dimensional surface is known as a \emph{slender body theory}. The first generation of methods in this category was called \emph{resistive force theory} \cite{hancock1953self}.  These methods were based on treating the fiber as a succession of prolate spheroids while ignoring nonlocal hydrodynamic interactions; that is, according to resistive force theory, the centerline velocity at a given location depends on the local force applied there, but not on the forces applied elsewhere on the fiber. Subsequent improvements accounted for nonlocal hydrodynamic effects \cite{cox1970motion,batchelor1970slender,lighthill1976flagellar}; a prominent example is the Keller-Rubinow formulation \cite{keller1976swimming}, reformulated for a periodic fiber \cite{shelley2000stokesian,cortez2012slender}: 
\begin{equation}
8\pi\mu \bm u_C^{KR}(s) 
=
\left[\left(I-3\bm e\bm e^T\right) - 2\log\left(\frac{\pi\epsilon}{4}\right)\left(I+2\bm e\bm e^T\right)\right] \bm f(s) + 
\int_{0}^{1}\left\{\left(\frac{I}{|\bm r|} 
+ \frac{\bm r\bm r^T}{|\bm r|^3}\right)\bm f(t) - \frac{I+\bm e(s) \bm e(s)^T}{|\sin(\pi(s-t))|/\pi}\bm f(s)\right\}
\,dt .
\label{eq:sbt}
\end{equation}
This equation assumes a unit-length closed fiber with arclength parameterization $\bm\gamma(s)$; from this we define $\bm r =  \bm\gamma(s) - \bm\gamma(t)$ and $\bm e = {\bm\gamma}'(s)$. The fiber radius is $\epsilon$, and $\mu$ is the fluid viscosity. The first term on the right-hand side is a local contribution to the centerline velocity from the force $\bm f$ imposed there, while the integral term represents the nonlocal hydrodynamic interactions with the rest of the fiber surface. 
This formulation, along with refinements that address the case where the fiber has free ends, has been widely used  \cite{johnson1980improved,tornberg2004simulating,nazockdast2017fast,lindner2015elastic,lauga2009hydrodynamics,free_ends,gotz2000interactions,li2013sedimentation}. 
However, the resistive force theory and its nonlocal successors all assume that the fiber does not closely approach itself or other fibers. In Keller and Rubinow's derivation \cite{keller1976swimming} this assumption is expressed in their choice of an expression for the velocity in the inner region using the method of matched asymptotic expansions. More generally, it has not been clear exactly which fully three-dimensional problem any of the prior slender body theories are approximating, that is, they are not derived from a boundary value problem (BVP) with explicit boundary conditions at the fiber surface. 

More recently, some of us presented a well-posed three-dimensional BVP  \cite{mori2020theoretical} given only the one-dimensional centerline force density, which we summarize as follows. We assume that the length scales for the flow problem are small enough that the Stokes model is appropriate. Writing $\bm u$ for velocity, $p$ for pressure, and $\mu$ for viscosity, the PDE for the fluid domain is: 
\begin{equation}
\label{eq:stokespde}\bm 0 = -\bm\nabla p + \mu \nabla^2 \bm u, \qquad 
0 = \bm\nabla\cdot\bm u .
\end{equation}
This PDE must be augmented by boundary conditions. We assume that the fluid is infinite in extent and quiescent, that is, $\bm u(\bm x)\to \bm 0$ as $|\bm x| \to \infty$. Then the only boundary of the fluid domain is the surface of a single immersed fiber whose centerline is a closed loop parameterized by a $C^2$ function from the circle $\mathbb{S}^1$ into $\mathbb{R}^3$. We assume a small constant fiber radius $\epsilon>0$ so that the two-dimensional surface of the fiber has the parameterization 
\begin{equation}
\bm X(s,\theta) = \bm \gamma(s) + \epsilon \cos(\theta) \bm N(s) + \epsilon \sin\theta \bm B(s) .
\end{equation}  
We call $s$ the \emph{centerline coordinate} and $\theta$ the \emph{circumferential coordinate}. The Bishop frame vectors $\bm N$ and $\bm B$ have unit length and are perpendicular to each other and to the tangent vector $\bm\gamma'(s)$; 
see Appendix A for an explanation of why we use this frame instead of the simpler Frenet-Serret frame as well as the details of our numerical implementation. We can now state the boundary conditions for the PDE. Let $\bm c(s)$ denote a surface velocity profile that depends only on $s$ (not on $\theta$), let $\bm \sigma(s,\theta)$ be the surface stress, let $\bm \nu$ denote the surface normal vector, and let $\bm f(s)$ be the imposed centerline force density. Then we require:
\begin{align}
\label{eq:fiberintegrity}\bm u(\bm X(s,\theta)) &= \bm c(s) , \\
\label{eq:avgforcebc}\int_{0}^{2\pi} \sigma(s,\theta) \cdot \bm \nu(s,\theta) J(s,\theta)\,d\theta &= \bm f(s) |\gamma'(s)| . 
\end{align}
The first equation \eqref{eq:fiberintegrity} restricts the boundary values of $\bm u$ by requiring that the fluid velocity be constant on each cross section of the fiber, a \emph{fiber integrity} constraint that preserves the circular shape of the cross sections. The centerline velocity function $\bm c(s)$ is unknown and must be determined as part of the solution. The second condition \eqref{eq:avgforcebc}, on the surface derivatives of $\bm u$, states that the integral of surface traction around each cross section must equal the imposed centerline force density $\bm f(s)$, which is given as part of the problem. On the left side of \eqref{eq:avgforcebc}, the stress tensor $\sigma$ is defined by the limiting value of $\sigma_{ij} = -p \delta_{ij} + \mu(\partial u_i/\partial x_j + \partial u_j/\partial x_i)$ as we approach the boundary from within the fluid, and the factor $J(s,\theta)$ is the surface Jacobian, $J(s,\theta) = |\bm X_s \times \bm X_\theta|$. Note that integration of either side of \eqref{eq:avgforcebc} from $s=0$ to $s=2\pi$ yields the net force exerted on the fiber by the fluid; the factor $|\bm\gamma'(s)|$ appears on the right side to account for the possibility of using a nonconstant parameterization speed. 

The boundary value problem defined by \eqref{eq:stokespde}-\eqref{eq:fiberintegrity}-\eqref{eq:avgforcebc}, first presented in \cite{mori2020theoretical}, presents some interesting challenges from a numerical perspective. This paper presents a new computational procedure appropriate for this new version of slender body theory. The new formulation is based on a representation of the velocity as the sum of a single-layer potential on the fiber surface and a distribution of point source singularities on its centerline; we call this the \emph{completed single-layer potential}. This formulation allows for convenient traction formulas while also circumventing the traditional disadvantages of the unmodified single-layer potential, as we discuss in Section \ref{sec:cslp}. We then present our discretization of the problem, wherein the unknowns are a finite number of Fourier coefficients of the single layer density, together with evenly spaced values of the centerline velocity $\bm c(s)$. The details of this procedure are covered in Section \ref{sec:discretization}, while the related matter of quadrature for singular integrals on the surface of a thin fiber appears in Appendix \ref{sec:quadratureappendix}. We study the convergence rate of our numerical method in Section \ref{sec:validation}, including a comparison to exact solutions for a uniformly translating rigid torus. These `exact' solutions were originally reported only to four-digit accuracy and so we also recomputed them to 12-digit or higher accuracy following the method of \cite{amarakoon1982drag} from 1982. Finally we compare our numerical method to the widely used and computationally efficient slender-body theory of Keller and Rubinow. 
The two methods agree to approximately order $\mathcal{O}(\epsilon^{1.7})$ when the centerline does not approach itself, a somewhat greater convergence rate than was rigorously proved in \cite{mori2020theoretical}, which demonstrated agreement at order $\mathcal{O}\left(\epsilon |\log \epsilon|^{3/2}\right)$. 
This $\mathcal{O}(\epsilon^{1.7})$ agreement is close to the $\mathcal{O}(\epsilon^2|\log\epsilon|)$ error predicted by the matched asymptotics of \cite{gotz2000interactions,keller1976swimming, johnson1980improved,tornberg2004simulating}.
We also present tests showing that the Keller-Rubinow slender body theory breaks down when the fiber surface comes close to self-intersection (Section \ref{sec:comp2sbt}). Thus, our method may be more suitable for simulating densely packed fiber suspensions. We conclude with priorities for future extensions of the method, such as treatments of multiple fibers, free ends, dynamic problems, and inertial flows. 

\section{The completed single-layer potential}
\label{sec:cslp}
Two attractive and widespread representations of unbounded Stokes flows are the single- and double- layer potentials: 
\begin{align}
\label{eq:slp} u_i(\bm x) &= \frac{1}{8\pi} \int_D \left(\frac{\delta_{ij}}{r} + \frac{r_ir_j}{r^3}\right) \rho_j(\bm y) \,dS_y , \\
\label{eq:dlp} u_i(\bm x) &= \frac{-3}{4\pi} \int_D \frac{r_ir_jr_k}{r^5}\nu_k(\bm y)\mu_j(\bm y) \,dS_y .
\end{align}
Here $\bm x$ is an observation point in the fluid domain, $D$ is the surface of an immersed particle, $\bm r = \bm x-\bm y$, the surface normal vector $\bm\nu$ points out of the particle into the fluid, and $\bm\rho$ and $\bm \mu$ are the single- and double-layer densities, respectively. We use Einstein's implicit summation notation (sum over the repeated indices $j$ and $k$). The tensor $\delta_{ij}/r + r_ir_j/r^3$ in \eqref{eq:slp} is the \emph{point force} or \emph{Stokeslet} on $\mathbb{R}^3$ so the single-layer potential is a sum of point forces located on the fiber surface. Similarly, $-6r_ir_jr_k/r^5$ is the \emph{stresslet} and the double-layer potential \eqref{eq:dlp} is a sum of point stresses located on the fiber surface. 
An arbitrary Stokes flow exterior to $D$ and decaying at infinity can be represented by a sum of these potentials; in this case the densities $\bm \rho$ and $\bm \mu$ have physical meaning (the surface traction and surface velocity, respectively). When only one or the other of these potentials appears, there are some flows which cannot be uniquely represented. For example, the single-layer potential is unable to represent flows with a nonzero volume flux across $D$, and if the density $\bm \rho$ is taken as a multiple of the normal vector $\bm \nu$ the resulting single-layer velocity field is zero. In the language of linear algebra, the single-layer potential as an operator mapping $\bm \rho$ to $\bm u$ has a one-dimensional nullspace, and its range has codimension one within the space of all quiescent Stokes flows evaluated at the fiber surface. The double-layer potential is also deficient: it cannot represent any flow which exerts a net force or torque on a closed particle, and the velocity in the fluid domain is zero if the density function $\bm \mu$ matches a rigid-body motion. The double-layer potential therefore has a six-dimensional nullspace and its range has codimension six.  Thus the single-layer and the double-layer potential, used in isolation, are each inadequate to represent general particulate flows. 

\begin{figure}
\includegraphics[height=2.2in]{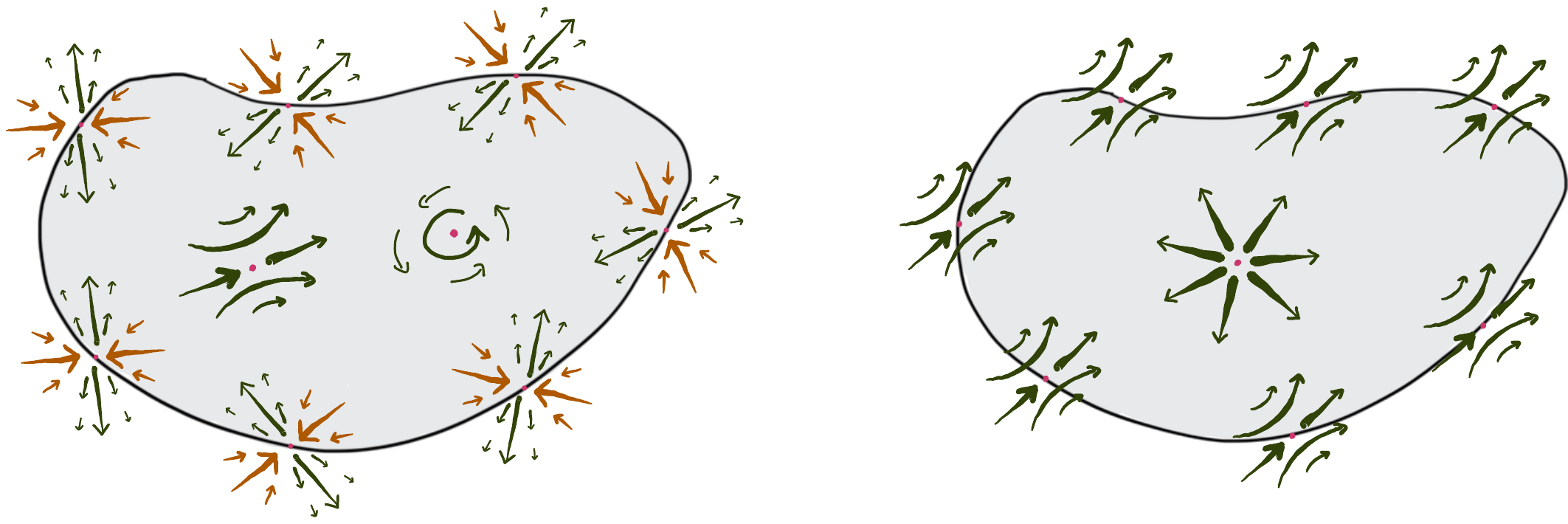}
\caption{ The completed single- and double-layer potential representations of Stokes flow are constructed by superposing fundamental solutions with singularities in the interior or distributed over the surface of an immersed particle. In Power and Miranda's completion of the double layer potential (at left), the flow consists of a distribution of stresslets on the surface together with a point force and a point torque in the interior (shaded). Our proposed completion of the single-layer potential appears on the right; the flow consists of a distribution of point forces on the surface together with one or more point sources in the interior (shaded).  The total strength of the internal point sources is proportional to the inner product of the surface normal with the single-layer density function.  
}
\label{fig:sldl}
\end{figure}

Power and Miranda provided a completion of the double-layer potential in the form of a point force and a point torque inside of an immersed solid object \cite{power1987second}. We sketch this formulation on the left side of Fig. \ref{fig:sldl}. The strengths of these internal singularities can be defined as six independent integrals of the double-layer density $\bm\mu$. This modified flow representation operator is now of full rank: any surface velocity field, including one which exerts a force and torque on the enclosed volume, can be uniquely represented. This procedure leads to a widely used Fredholm integral equation of the second kind for solving Dirichlet problems. The completed double-layer potential does have a disadvantage despite its wide use: if one needs the local surface traction (rather than the net force), it generally requires evaluating hypersingular integrals of the density, although there are remedies for this problem in the case of rigid particle motions \cite{keaveny2011applying,corona2017integral,mitchell2017generalized}. 

The single-layer potential also has a long computational history. The first numerical implementation of a boundary integral method was based on an unmodified single-layer representation; see \cite{youngren_acrivos_1975}. 
The single-layer representation has the advantage that the local surface tractions are comprised of convergent integrals; however, there are two problems with the unmodified single-layer formulation. 
The first is the existence of a nullspace (a zero eigenvalue) discussed above. The second is that the operator carrying the density $\bm \rho$ to the surface velocity $\bm u$ has arbitrarily small nonzero eigenvalues. In a numerical implementation, any surface integration procedure is subject to numerical error, which can lead to invertible discrete linear operators even when the continuous operator is singular. The condition number of the discrete operator increases as the integration procedure becomes more accurate, but the method can give acceptable results in an intermediate range where the discretization is neither too coarse nor too fine. For problems where the surface velocity is prescribed, the completed double-layer methods lead to second-kind integral equations, which do not suffer from this conditioning issue. 

In the more specific setting where the particle is a slender fiber, Koens and Lauga recently obtained versions of slender-body theory by starting with either the single- or double-layer potentials and using matched asymptotic expansions in the fiber radius \cite{koens2018boundary}. They found that the use of a single-layer potential leads to a singular system of equations; in particular, the first Fourier mode of the force density in the circumferential direction is not uniquely determined from the surface velocity, while no solution exists if the first mode of the surface velocity is nonzero. As they remark, this failure corresponds to the inability of the single-layer potential to represent volume changes. 

We now offer a completion procedure for single-layer potentials which resolves the nullspace issue and may have some advantages in comparison to the double-layer formulations. The procedure is very simple: in addition to the single-layer potential, we include the flow due to one or more point sources located within the particle whose total strength is proportional to the inner product of the single-layer density with the surface normal vector. More concretely, let $D$ be a closed particle surface, and let $\bm \chi$ be a map carrying points on the surface to points in the interior. For a convex particle, $\bm \chi$ could be a constant map carrying surface points to the particle centroid. In our setting with a closed slender fiber, we let $\bm\chi$ carry surface points to the corresponding centerline point. We then define the velocity on the exterior to the particle by the modified single-layer equation 
\begin{equation}
u_i(\bm x) = \frac{1}{8\pi}\int_D \left( \frac{\delta_{ij}}{|\bm x - \bm y|} + \frac{(x_i-y_i)(x_j-y_j)}{|\bm x - \bm y|^3}\right) \rho_j(\bm y) \,dS_{\bm y} + \frac{1}{4\pi}\int_D \frac{\nu_j(\bm y)\rho_j(\bm y)(x_i-\chi_i(\bm y))}{|\bm x - \bm \chi(\bm y)|^3}\,dS_{\bm y} .
\label{eq:modsing}
\end{equation}
The first term in this equation is the single-layer potential with density $\bm \rho$. The second term is the sum of the velocities at $\bm x$ due to point sources distributed at $\bm \chi(\bm y)$, each with strength $\bm \nu(\bm y) \cdot \bm \rho(\bm y)$, where $\bm \nu$ is the normal vector pointing out of the particle (into the fluid domain). We note that the first integrand becomes singular if $\bm x$ is on the particle surface, but the second integral remains regular even in that case. We are using the Green's functions for unbounded flow, although the idea should also work for bounded flow domains. Although this procedure is exactly analogous to the widely used completion procedure for the double layer formulations, we have not found any discussion of it in the literature.

The main advantage of the completed single-layer representation is that the surface tractions are easy to compute. The single-layer potential has a known stress field, and so does the point source. At observation points $\bm x$ in the bulk fluid the stress tensor is 
\begin{equation}
\begin{split}
\sigma_{ik} & = \frac{-3}{4\pi}\int_D \frac{(x_i - y_i)(x_j-y_j)(x_k-y_k)}{|\bm x - \bm y|^5} \rho_j(\bm y) \,dS_y  \\
& \qquad +  \frac{1}{4\pi}\int_D \left(\frac{2\delta_{ik}}{|\bm x - \bm \chi(\bm y)|^3} - 6\frac{(x_i-\chi_i(\bm y))(x_k-\chi_k(\bm y))}{|\bm x - \bm \chi(\bm y)|^5} \right)\nu_j(\bm y)\rho_j(\bm y)   \,dS_{\bm y} ,
\end{split}
\end{equation}
and the surface traction is the contraction of the stress with the surface normal in the limit as $\bm x$ approaches the fiber surface from the fluid domain, that is, from the side into which $\bm\nu$ points. For the point source the integrands remain regular and for the single-layer part we use results from Pozrikidis \cite{pozrikidis1992boundary} to arrive at 
\begin{gather}
\begin{split}
t_i(\bm x) & = -\frac12\rho_i(\bm x) - \frac{3}{4\pi}\int_D \frac{(x_i-y_i)(x_j-y_j)(x_k-y_k)}{|\bm x-\bm y|^5} \nu_k(\bm x) \rho_j(\bm y)\,dS_{\bm y}\\
& \quad +  \frac{1}{4\pi}\int_D \left(\frac{2\delta_{ik}}{|\bm x - \bm \chi(\bm y)|^3} - 6\frac{(x_i-\chi_i(\bm y ))(x_k-\chi_k(\bm y ))}{|\bm x - \bm \chi(\bm y)|^5} \right)\rho_j(\bm y )\nu_k(\bm x)\nu_j(\bm y)   \,dS_{\bm y }.
\label{eq:tractioncslpgeneral}
\end{split}
\end{gather} 

When $\bm{\chi}$ is a constant function, we can prove that the velocity representation \eqref{eq:modsing} is unique and can represent an arbitrary flow. We expect but do not prove that it also holds for nonconstant $\bm\chi$, which is a more computationally appropriate choice for the slender fiber geometry. 
In the constant case the velocity can be more simply written as 
\begin{equation}
u_i(\bm x) = \frac{1}{8\pi}\int_D \left( \frac{\delta_{ij}}{|\bm x - \bm y|} + \frac{(x_i-y_i)(x_j-y_j)}{|\bm x - \bm y|^3}\right) f_j(\bm y ) \,dS_{\bm y } + \frac{1}{4\pi}\frac{(x_i-\chi_i(\bm y ))}{|\bm x - \bm \chi(\bm y)|^3}\int_D \nu_j(\bm y)f_j(\bm y )\,dS_{\bm y }
\end{equation}
so that the rate of volume creation at $\bm \chi(\bm y)$ is precisely the inner product of $\bm \nu$ and $\bm f$. Thus, to represent an arbitrary flow, one first determines the volume flux rate $\alpha$ and then takes an initial surface distribution $\bm f_0 = \frac{\alpha}{|D|}\bm \nu$. Then the velocity induced by $\bm f_0$ has the desired volume flux. Now the difference between this flow and the desired one is flux-free and can be represented in infinitely many ways by a single-layer potential, but only in one way by a single-layer potential with density $\bm f_1$ satisfying $\int_D \bm f_1 \cdot \bm \nu = 0$.  Then $\bm f = \bm f_0 + \bm f_1$ is a density function inducing the desired velocity. Uniqueness is also a consequence of the fact that a flow with zero flux can be uniquely represented by a single-layer potential whose density has zero inner product with the surface normal.  


In the remainder of this paper we use this completed single-layer velocity formulation to develop a computational method suitable for simulating closed slender fibers. 
In a resonance with the finding by Koens and 
Lauga that the slender body theory based on a single-layer potential has a deficiency at mode one \cite{koens2018boundary}, we find that our correction procedure modifies only the terms corresponding to modes zero, one, and two in the discrete version, and the greatest modification is to the first mode.

\section{Discretization of the slender body BVP}
\label{sec:discretization}
\begin{figure}
\[\includegraphics[width=0.4\linewidth]{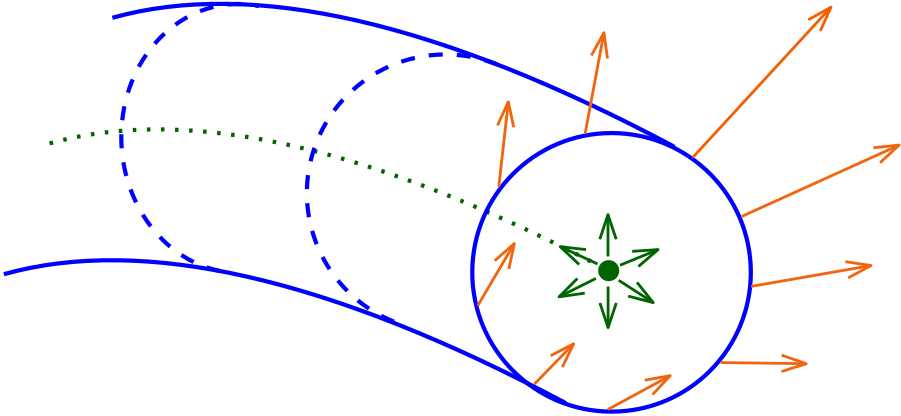}\]
\caption{We represent the fluid velocity using a distribution of point forces on the fiber surface (orange arrows) together with a distribution of point sources on the the centerline (green arrows). The strength of the point sources are taken equal to the fluxes of the single layer density through each cross section, so that the volume creation rate at $\bm\gamma(s)$ is $\int_0^{2\pi} \bm \nu(s,\theta) \cdot \bm \rho(s,\theta) J(s,\theta)\,d\theta$. In the figure, the orange arrows have a positive flux out of the cross section and accordingly the green arrows indicate a point source rather than a sink.} 
\label{fig:cross_section_sketch}
\end{figure}

As stated previously, we consider a single fiber in a quiescent fluid without boundary. The fiber centerline is a closed loop. The fluid velocity is represented by the sum of a single-layer potential and a distribution of point sources along the centerline; see Fig. \ref{fig:cross_section_sketch}. 
Let $\bm \gamma:[0,2\pi)\to\mathbb{R}^3$ be a parameterization of the centerline, not necessarily of constant speed. 
Let $\bm X(s,\theta)$ parameterize the fiber surface. The circular cross sections are normal to the centerline $\bm \gamma$ and have uniform radius $\epsilon$. 
Let $\bm \nu(s,\theta)$ denote the unit surface normal vector pointing out of the fiber (into the fluid domain), and let $J(s,\theta)$ denote the Jacobian, $J = |\bm X_s\times \bm X_\theta|$. 
Then the completed single-layer fluid velocity equation \eqref{eq:modsing} reduces to:
\begin{equation}u_i(\bm x) = \frac{1}{8\pi}\int_D \left(\frac{\delta_{ij}}{r}+\frac{r_ir_j}{r^3} \right)\rho_j(\bm y)\,dS_y 
\label{eq:SLP}
+ \frac{1}{4\pi}\int_0^{2\pi} \frac{R_i}{R^3}\int_0^{2\pi} \rho_j(s,\theta)\nu_j(s,\theta)J(s,\theta)\,d\theta \,ds ,
\end{equation}
where the vector $\bm r = \bm x-\bm y$ points from the surface integration point to the observation point, and the vector $\bm R = \bm x - \bm \gamma(s)$ points from the centerline integration point to the observation point. The first term is a single-layer potential with density $\bm \rho$, and the second term is a distribution of point sources over the centerline.

Similarly, the equation \eqref{eq:tractioncslpgeneral} for the surface traction exerted by the exterior fluid on the fiber at a surface point $\bm x=\bm X(s^*,\theta^*)\in D$ becomes 
\begin{equation}
\begin{split}
t_i(\bm x) =& -\frac12 \rho_i(\bm x) - \frac{3}{4\pi}\int_0^{2\pi} \int_0^{2\pi}
\frac{r_ir_jr_k}{r^5} \rho_j(\bm y)\nu_k(\bm x)J(s,\theta)\,d\theta\,ds \\
&+ \frac{1}{4\pi} \int_0^{2\pi} \left(\frac{2\nu_i(s^*,\theta^*)}{R^3} - 6 R_i\frac{R_j \nu_j(s^*,\theta^*)}{R^5}\right)\int_0^{2\pi} \rho_k(s,\theta)\nu_k(s,\theta)J(s,\theta)\,d\theta\,ds
\end{split}
\end{equation}
with $\bm r$ and $\bm R$ now given by 
$\bm r = \bm X(s^*,\theta^*) - \bm X(s,\theta)$ and $\bm R = \bm X(s^*,\theta^*) - \bm \gamma(s)$. 

We can now substitute these velocity and traction expressions into the boundary conditions \eqref{eq:fiberintegrity}-\eqref{eq:avgforcebc}.  
Let $\bm f(s)$ be a vector function of the centerline describing the force density on a cross section. Then we have 
\begin{align}
0 =& -c_i(s^*)+\frac{1}{8\pi} \int_0^{2\pi}\int_0^{2\pi} 
\left(\frac{\delta_{ij}}{r}+\frac{r_ir_j}{r^3}\right)\rho_j(s,\theta)J(s,\theta)
\,ds\,d\theta + \frac{1}{4\pi}\int_0^{2\pi} \frac{R_i}{R^3}\int_0^{2\pi} \rho_j(s,\theta)\nu_j(s,\theta)J(s,\theta)\,d\theta \,ds
\label{eq:velfunarc}
\\
\begin{split}
f_i(s^*) \left|\bm \gamma'(s^*)\right|=&  \int_0^{2\pi} 
\left[
-\frac12 \rho_i(s^*,\theta^*) - 
\frac{3}{4\pi} \int_0^{2\pi}\int_0^{2\pi} 
\frac{r_ir_jr_k}{r^5}\rho_j(s,\theta)\nu_k(s^*,\theta^*)J(s,\theta)
\,ds\,d\theta\right.\\
&+\left. \frac{1}{4\pi} \int_0^{2\pi} \left(\frac{2\nu_i(s^*,\theta^*)}{R^3} - 6 R_i\frac{R_j \nu_j(s^*,\theta^*)}{R^5}\right)\int_0^{2\pi} \rho_k(s,\theta)\nu_k(s,\theta)J(s,\theta)\,d\theta\,ds\right] J(s^*,\theta^*)\,d\theta^*,
\label{eq:avgforce}
\end{split}
\end{align}
where the unknowns are the velocity $\bm c(s)$, a function of the centerline only, and the single layer density $\bm \rho(s,\theta)$. 
The first equation says that the surface velocity is independent of the circumferential coordinate $\theta$ (the \emph{fiber integrity condition}). The second condition states that the integral of the surface traction around a circular cross section of the fiber matches $\bm f$. The factor of $\left|\bm \gamma'(s^*)\right|$ on the left of \eqref{eq:avgforce} accounts for the speed of the parameterization so that upon integration from $0$ to $2\pi$ in $s^*$, we would find that the net force (right side) is equal to the integral of the centerline force with respect to arclength (left side).

To discretize this system of integral equations, we start by enforcing them only at finitely many points. 
For \eqref{eq:velfunarc} we let $(s^*,\theta^*)$ range over a regular grid. 
Letting $n_s$ and $n_\theta$ be odd integers giving the number of grid points in each direction, we take $(s^*,\theta^*) = (2\pi j_s/n_s,2\pi j_\theta/n_\theta)$ for $0\le j_s < n_s$ and $0\le j_\theta < n_\theta$.   
Similarly, we enforce \eqref{eq:avgforce} only for $s^* = 2\pi j_s/n_s$ with $0\le j_s < n_s$. 
This gives a total of $3(n_sn_\theta + n_s)$ scalar equations. 
To obtain finitely many unknowns, we seek the density $\bm \rho$ in a finite-dimensional space of complex exponentials:
\[\rho_\ell(s,\theta) = \sum_{k_s,k_\theta} \alpha_{\ell, k_s, k_\theta} 
\exp\left(\sqrt{-1}(k_s s + k_\theta \theta )\right) . 
\]
Here the indices range over $ - (n_s-1)/2 \le k_s \le (n_s-1)/2  $ and $ -(n_\theta-1)/2  \le k_\theta \le (n_\theta-1)/2$. Considering the three possible values of the space dimension index $\ell$, we have a total of $3n_sn_\theta$ unknown Fourier coefficients of $\bm \rho$. The values of $c_i(s^*)$ at $s^* = 2\pi j_s/n_s$ provide the remaining $3 n_s$ unknowns, which we abbreviate by writing $c_{i,j_s} = c_i\left(\frac{2\pi j_s}{n_s}\right)$ for $0\le j_s< n_s$. Upon substituting this expression for $\bm \rho$ and moving the sums outside the integrals, we obtain the following discrete equations: 
\begin{align}
\label{eq:dv}
0=-c_{i,j_s}+ \frac{1}{8\pi}&\sum_{j,k_s,k_\theta} \alpha_{j,k_s,k_\theta} \int_0^{2\pi}\int_0^{2\pi} \left(\frac{\delta_{ij}}{r} + \frac{r_ir_j}{r^3} + 2\frac{R_i\nu_j(s,\theta)}{R^3}\right)\exp\left(\sqrt{-1}(k_ss+k_\theta\theta)\right)J(s,\theta)\,d\theta\,ds\\
F_{i,j_s}  =\nonumber \frac{2\pi}{n_\theta} \sum_{j_\theta =0}^{n_\theta-1}J(s^*,\theta^*) &
\Biggl\{
\sum_{j,k_s,k_\theta} \alpha_{j,k_s,k_\theta}
\Biggl[
-\frac12\delta_{ij}\exp\left(\sqrt{-1}\left( k_ss^*+k_\theta \theta^* \right)\right) 
+ \frac{1}{4\pi}\int_0^{2\pi}\int_0^{2\pi}
\Biggl(-3\frac{r_ir_j(\bm r\cdot\bm \nu(s^*,\theta^*))}{r^5}\\
&+ 2\frac{\nu_i(s^*,\theta^*)}{R^3}\nu_j(s,\theta) - 6R_i\nu_j(s,\theta)\frac{\bm R \cdot \bm \nu(s^*,\theta^*)}{R^5}
\Biggr)\exp\left(\sqrt{-1}(k_ss+k_\theta\theta)\right)
J(s,\theta)\,ds\,d\theta
\Biggr]
 \Biggr\} . \label{eq:df}
\end{align}
In these equations we are using the abbreviations $s^* = 2\pi j_s/n_s$, $\theta^* = 2\pi j_\theta/n_\theta$, $F_{i,j_s} = f_i(s^*)|\bm\gamma'(s^*)|$, $\bm r = \bm X(s^*,\theta^*) - \bm X(s,\theta)$, and $\bm R = \bm X(s^*,\theta^*) - \gamma(s)$.
The discrete fiber integrity equation \eqref{eq:dv} holds for all $0\le j_s < n_s$ and $0\le j_\theta < n_\theta$, while the discrete averaged force equation \eqref{eq:df} holds for $0\le j_s\le n_s $. Note that to arrive at \eqref{eq:df} we have replaced the outermost integral (in $\theta^*$) from \eqref{eq:avgforce} with a trapezoidal rule sum. 

To set up the linear algebra system, we have to evaluate the integrals in \eqref{eq:dv}-\eqref{eq:df}. 
The Stokeslet integrand in \eqref{eq:dv} and the stresslet integrand in \eqref{eq:df} both have a $1/r$ singularity as $(s,\theta) \to (s^*,\theta^*)$, so the numerical quadrature procedure is a nontrivial problem. We give details of our method in Appendix \ref{sec:quadratureappendix}. An interesting feature of this formulation is that the accuracy of the quadrature can be chosen independently of the matrix size. Once the integrals have been computed, we have a dense and non-normal system of linear equations. The condition numbers for the problems we considered range from approximately $10^3$ to $10^8$. 
We chose to use the SVD for the linear solve because of its good performance with poorly conditioned systems; the computational expense of this method for our dense, non-normal system of linear equations is acceptable because the overall solution time is dominated by the matrix assembly rather than the linear solve.

\begin{table}\centering
\ra{1.1}
\begin{tabular}{@{}llccc@{}}\toprule
&& $n_\theta=7$ & $n_\theta=13$ & $n_\theta=25$ \\ \midrule
$\epsilon=10^{-2}$ &
$n_s = 7$ & $1.2026\cdot 10^3$ & $2.4031\cdot 10^3$ & $4.8038 \cdot 10^3$\\
&$n_s = 21$ & $1.5546 \cdot 10^3$ & $2.4139 \cdot 10^3$ & $4.8190 \cdot 10^3$\\ 
&$n_s = 63$ & $1.3127 \cdot 10^5$ & $1.2857 \cdot 10^5$ & $1.2705 \cdot 10^5$\\ 
&$n_s = 189$ & $1.8706 \cdot 10^4$ & $1.6422 \cdot 10^4$ & $1.5006 \cdot 10^4$\\ \midrule
$\epsilon=10^{-3}$ &
$n_s = 7$ & $1.2002\cdot 10^4$ & $2.4002\cdot 10^4$ & $4.8002 \cdot 10^4$\\
&$n_s = 21$ & $1.3306 \cdot 10^4$ & $2.4007 \cdot 10^4$ & $4.8002 \cdot 10^4$\\ 
&$n_s = 63$ & $4.8891 \cdot 10^4$ & $4.8882 \cdot 10^4$ & $4.8008 \cdot 10^4$\\ 
&$n_s = 189$ & $2.7247 \cdot 10^5$ & $2.7180 \cdot 10^5$ & $2.7142 \cdot 10^5$\\ \midrule
$\epsilon=10^{-4}$ &
$n_s = 7$ & $1.2000 \cdot 10^5$ & $2.4000\cdot 10^5$ & $4.8000 \cdot 10^5$\\
&$n_s = 21$ & $1.2991 \cdot 10^5$ & $2.4001 \cdot 10^5$ & $4.8001 \cdot 10^5$\\ 
&$n_s = 63$ & $3.8928 \cdot 10^5$ & $3.8930 \cdot 10^5$ & $4.8002 \cdot 10^5$\\ 
&$n_s = 189$ & $1.5466 \cdot 10^6$ & $1.5466 \cdot 10^6$ & $1.5466 \cdot 10^6$\\ 
 \bottomrule
\end{tabular}
\caption{The condition numbers of the discrete linear systems increase as the fiber radius $\epsilon$ decreases, and they do not increase as the quadrature is refined. The condition numbers generally increase with the centerline and circumferential discretization parameters $n_s$ and $n_\theta$, with the exception of the anomalous row corresponding to $n_s=189$, $\epsilon = 10^{-2}$. The discrete systems reported in this table were generated for a trefoil knot centerline using a well resolved quadrature ($q_n=40$ for all reported values; see Appendix \ref{sec:quadratureappendix} for details). Further tests with a more refined quadrature $q_n=50$ (not reported here) gave relative changes of less than $1/10000$ compared to those appearing in this table. 
}
\label{tbl:conds}
\end{table}

\subsection{Circumferential integrals of nonsingular terms} 
The variable $R = |\bm X(s^*,\theta^*) - \gamma(s)|$ has no dependence on the circumferential integration variable $\theta$; moreover, its minimum as a function of $s$ is $\epsilon$ rather than zero. Therefore the terms with $R$ in the denominators are more analytically tractable than those involving $r = |\bm X(s^*,\theta^*) - \bm X(s,\theta)|$. In particular, they can always be reduced to one-dimensional integrals and for many parameter values they simply vanish. To see this, define the circumferential integrals
\begin{equation}
\mathcal{M}(k_\theta,s) = \int_0^{2\pi} \nu_j(s,\theta)\exp\left(\sqrt{-1}(k_ss+k_\theta\theta)\right)J(s,\theta)\,d\theta
\end{equation}
and use this notation to rewrite the $R$ integrals as 
\begin{gather}
\label{eq:innersl}\frac{1}{8\pi}\int_0^{2\pi}\int_0^{2\pi} 2\frac{R_i\nu_j(s,\theta)}{R^3}\exp\left(\sqrt{-1}(k_ss+k_\theta\theta)\right)J(s,\theta)\,d\theta\,ds 
=\frac{1}{4\pi} \int_0^{2\pi}  \frac{R_i}{R^3} \mathcal{M}(s,k_\theta)\,ds
\\
\begin{split}
\label{eq:innerdl}\frac{1}{4\pi}\int_0^{2\pi}\int_0^{2\pi}
\Biggl( 2\frac{\nu_i(s^*,\theta^*)}{R^3}\nu_j(s,\theta) - 6R_i\nu_j(s,\theta)\frac{\bm R \cdot \bm \nu(s^*,\theta^*)}{R^5}
\Biggr)\exp\left(\sqrt{-1}(k_ss+k_\theta\theta)\right)
J(s,\theta)\,ds\,d\theta\\
=\frac{1}{2\pi} \int_0^{2\pi} \left(\frac{\nu_i(s^*,\theta^*)}{R^3} - 3R_i\frac{\bm R \cdot \bm \nu(s^*,\theta^*)}{R^5}\right) \mathcal{M}(s,k_\theta)\,ds .
\end{split}
\end{gather}
Now the inner integrals can be solved explicitly by expanding 
$J(s,\theta) = \epsilon |\gamma'(s)| (1-\epsilon\cos\theta\kappa_1(s)-\epsilon\sin\theta\kappa_2(s))$ and 
$\bm \nu(s,\theta) = \cos\theta\bm N(s) +\sin\theta\bm B(s)$; this leads to \begin{equation}
\mathcal{M}(k_\theta,s)
= \begin{cases}
-\pi \epsilon^2 |{\bm\gamma}'(s)| \Big(N_j(s)\kappa_1 + B_j(s)\kappa_2\Big)e^{\sqrt{-1}k_s s} & k_\theta = 0\\
\pi\epsilon |{\bm\gamma}'(s)| \Big(N_j(s) + \sqrt{-1}B_j(s)\Big)e^{\sqrt{-1}k_s s} & k_\theta = 1\\
-\frac12 \pi\epsilon^2|{\bm\gamma}'(s)|\Big(N_j(s) + \sqrt{-1}B_j(s)\Big) \Big(\kappa_1+\sqrt{-1}\kappa_2\Big)e^{\sqrt{-1}k_s s}& k_\theta = 2\\
0&k_\theta>2 .
\end{cases} 
\label{eq:Mkths}
\end{equation}
The values for negative $k_\theta$ are the complex conjugates of these. Thus, we could use one-dimensional quadrature for $|k_\theta|<3$ and omit these terms entirely for $|k_\theta|\ge 3.$ The integrals which reduce so conveniently to one-dimensional quadrature problems all correspond to the point sources we introduced to complement the single-layer potential. As we mentioned in the previous section, a recent study \cite{koens2018boundary} shows that a slender-body formulation based on an unmodified single-layer potential will be noninvertible precisely at mode $k_\theta=1$.  It is encouraging to see that our correction procedure modifies the first mode most (order $\epsilon$) while leaving the $|k_\theta|>2$ terms unchanged and modifying the zeroth and second modes only by a factor of $\epsilon^2$. 
%

The reduction to one-dimensional quadrature for $|k_\theta|\le 2$ does not improve the speed of the overall algorithm, because evaluations on a two-dimensional quadrature grid are still needed for the singular integrals. In the current implementation, we use 2D quadrature for all integrals, but we omit the $R$ terms when $|k_\theta|> 2$, in accordance with \eqref{eq:Mkths}.

\section{Error analysis of the numerical method}
\label{sec:validation}
\subsection{Vertically translating torus}
In order to compare our method to an exact solution, we consider a torus whose centerline is a unit circle and whose cross sections have radius $\epsilon$, as above.  If the forcing is uniform and aligned with the rotational symmetry axis of the body, the resulting velocity is also uniform (it does not vary along the centerline) and in the same direction. Because of the simple geometry and because the computed velocity is a rigid-body motion, we can compare to an analytical solution due to Amarakoon \cite{amarakoon1982drag} and predecessors \cite{tchen1954motion,goren1980asymmetric,majumdar1977axisymmetric} which comes from separating variables in toroidal coordinates. As in that work, we investigate the behavior of a nondimensionalized force as $\epsilon$ varies. The quantity we evaluate is
\[F'_\infty = F'_\infty(\epsilon) = \frac{F}{6\pi\mu U (1+\epsilon)}.\] 
Here $U$ is the velocity of the torus and $F$ is the net force, both measured in the direction of the symmetry axis. The denominator is the force on a sphere whose outer diameter $(2+2\epsilon)$ coincides with the outer diameter of the torus surface, so we expect $F'_\infty<1$ and $F'_\infty\to 0$ as
 $\epsilon\to 0$. To be clear, the analytical solution \cite{amarakoon1982drag} is posed as a resistance problem (set $U=1$ and solve for $F$) while our method is essentially a mobility problem (set $F = 1$ and compute $U$), but the nondimensionalization allows us to compare the results. 

\begin{table}
\begin{tabular}{llll}
\toprule
$1/\epsilon$ &1982 result & Our exact value & Accuracy estimate\\
\midrule
$1.5$ &$0.9199  $& $0.9199051996705850$ & $5.01\cdot10^{-12}$ \\
$2.0$ &$0.9071  $& $0.9071647049517928$ & $4.36\cdot10^{-12}$ \\
$3.0$ &$0.88456  $& $0.8845567193948104$ & $9.95\cdot10^{-13}$ \\
$4.0$ &$0.86465  $& $0.8646521851778260$ & $5.59\cdot10^{-12}$ \\
$6.0$ &$0.83163  $& $0.8316248598100708$ & $8.07\cdot10^{-13}$ \\
$8.0$ &$0.80552  $& $0.8055222447767347$ & $2.40\cdot10^{-12}$ \\
$10.0$ &$0.78431  $& $0.7843079118776118$ & $1.68\cdot10^{-13}$ \\
$15.0$ &$0.74480  $& $0.7447973723171650$ & $7.20\cdot10^{-13}$ \\
$20.0$ &$0.71680  $& $0.7168018001135357$ & $4.14\cdot10^{-14}$ \\
$30.0$ &$0.67839  $& $0.6783917926050503$ & $1.43\cdot10^{-12}$ \\
$40.0$ &$0.65230  $& $0.6522959699851174$ & $1.44\cdot10^{-13}$ \\
$60.0$ &$0.61752  $& $0.6175207328846227$ & $5.56\cdot10^{-15}$ \\
$80.0$ &$0.59438  $& $0.5943759220857434$ & $8.38\cdot10^{-16}$ \\
$100.0$ &$0.57731  $& $0.5773112010340116$ & $1.78\cdot10^{-12}$ \\
$150.0$ &$0.54825  $& $0.5482504902513438$ & $1.56\cdot10^{-13}$ \\
$200.0$ &$0.52909  $& $0.5290943490801049$ & $2.78\cdot10^{-14}$ \\
$300.0$ &$0.50402  $& $0.5040229530606322$ & $2.16\cdot10^{-15}$ \\
\midrule
$10.0$ &$0.78431$& $0.7843079118776118$ & $1.68\cdot10^{-13}$ \\
$100.0$ &$0.57731$& $0.5773112010340116$ & $1.78\cdot10^{-12}$ \\
$1000.0$ && $0.4410799504734741$ & $0.0$ \\
$10000.0$ && $0.3552543625215476$ & $5.65\cdot10^{-16}$ \\
$100000.0$ && $0.2972352562536550$ & $4.71\cdot10^{-17}$ \\
\bottomrule
\end{tabular}
\caption{ Drag coefficients for a torus of centerline radius $1$ and cross-sectional radius $\epsilon$ translating along its symmetry axis; reported values are the net force scaled on $6\pi\mu U (1+\epsilon)$. Previously reported values were obtained in the 1980s and only listed five digits \cite{amarakoon1982drag}, so we have recomputed them. {
The accuracy estimates on the right refer to errors in these `exact' solutions, not comparisons to our numerical method.
} 
}
\label{tbl:exact1982}
\end{table}

Before discussing the convergence rate of our numerical method, we will digress to comment further on the reference solutions. Although we refer to the results published in the 1982 paper \cite{amarakoon1982drag} as ``exact,'' they were only reported to four or five digit accuracy. This is because the exact solution procedure yields an infinite set of equations relating the coefficients of the toroidal harmonic expansion of the solution, and this system of equations has to be truncated and then solved with a numerical linear algebra procedure. Thus even the ``exact'' solutions are subject to some numerical uncertainty. To get more digits, we reimplemented the procedure described by \cite{amarakoon1982drag} and used larger linear systems and double- rather than single-precision arithmetic. 
To assess the numerical error in these ``exact'' solutions, we computed the difference between the two sides of the identity $\sum_{n\ge 1} n B_n = \sum_{n\ge 0} C_n$, which appears on the first page of \cite{amarakoon1982drag} and report this quantity as a proxy for the accuracy of the results. Our recomputations of the exact solutions are given in Table \ref{tbl:exact1982} together with these accuracy estimates. These results, obtained through a strict reimplementation of the 1982 presentation, are sufficient to get as many digits as we require for verification of our own numerical method. However, we note for completeness that O'Neill and Bhatt gave an improved version \cite{o1991slow} of Goldman, Cox and Brenner's classic work on the problem of a sphere moving near a plane wall \cite{goldman1967slow}, removing the need to solve a linear system of equations. It is likely that similar methods could improve Amarakoon's formulation \cite{amarakoon1982drag}. However, we found that the unmodified 1982 algorithm already gives sufficiently accurate results, and we did not attempt to improve it. 

\begin{figure}
\includegraphics[width=\linewidth]{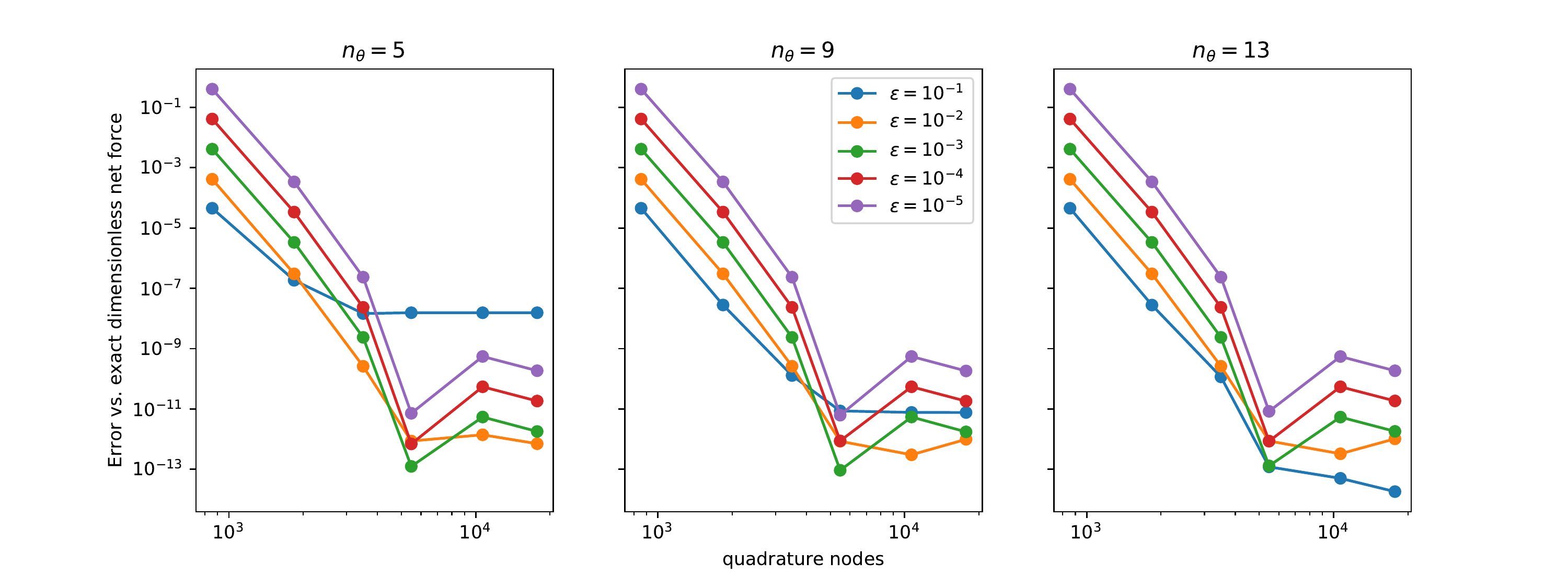}
\caption{A comparison of our numerical results to exact solutions for a simple torus with uniform forcing along the symmetry axis. The horizontal coordinate gives the number of quadrature nodes used to evaluate each matrix entry. The vertical axis shows the relative error versus the exact solutions from Table \ref{tbl:exact1982}. 
The method appears to be spectrally accurate with respect to the quadrature, subject to a ceiling on the accuracy imposed by the condition number (approximately $100/\epsilon$) and the grid resolution parameters $n_s$ and $n_\theta$. To obtain a single scheme we should have the quadrature rule and the parameters $n_s$ and $n_\theta$ increase together, but the optimal way to do this will depend on the fiber geometry and forcing. We note that a small $n_\theta$ is acceptable when $\epsilon$ becomes small.  
} 
\label{fig:vtt}
\end{figure}

In Fig. \ref{fig:vtt} we compare the accuracy of our numerical scheme to these recomputed exact solutions. For fixed grid parameters $(n_s, n_\theta)$, the convergence rate is spectral in the number of quadrature nodes used to evaluate the integrals, but it eventually reaches a plateau where other sources of error dominate. One source of error is from the presence of unresolved Fourier modes in the single-layer density $\bm \rho$; this is clearly the issue for the $\epsilon=10^{-1}$ curve in the left panel of Fig. \ref{fig:vtt}, because the issue is resolved by increasing the circumferential discretization parameter $n_\theta$.  For smaller values of the radius $\epsilon$, there is apparently no error reduction at all when $n_\theta$ increases from $5$ to $13$; the higher circumferential modes are, unsurprisingly, irrelevant, which allows us to save some computational expense in the more complicated simulations discussed later. When the fiber centerline is more complicated than a simple circle, the centerline discretization parameter $n_s$ will also play an important role. The condition number of our linear systems is approximately $100/\epsilon$, and so there is also a loss of accuracy as $\epsilon\to0$. 

\subsection{Trefoil knot} 
\begin{figure}
\includegraphics[width=0.40\linewidth]{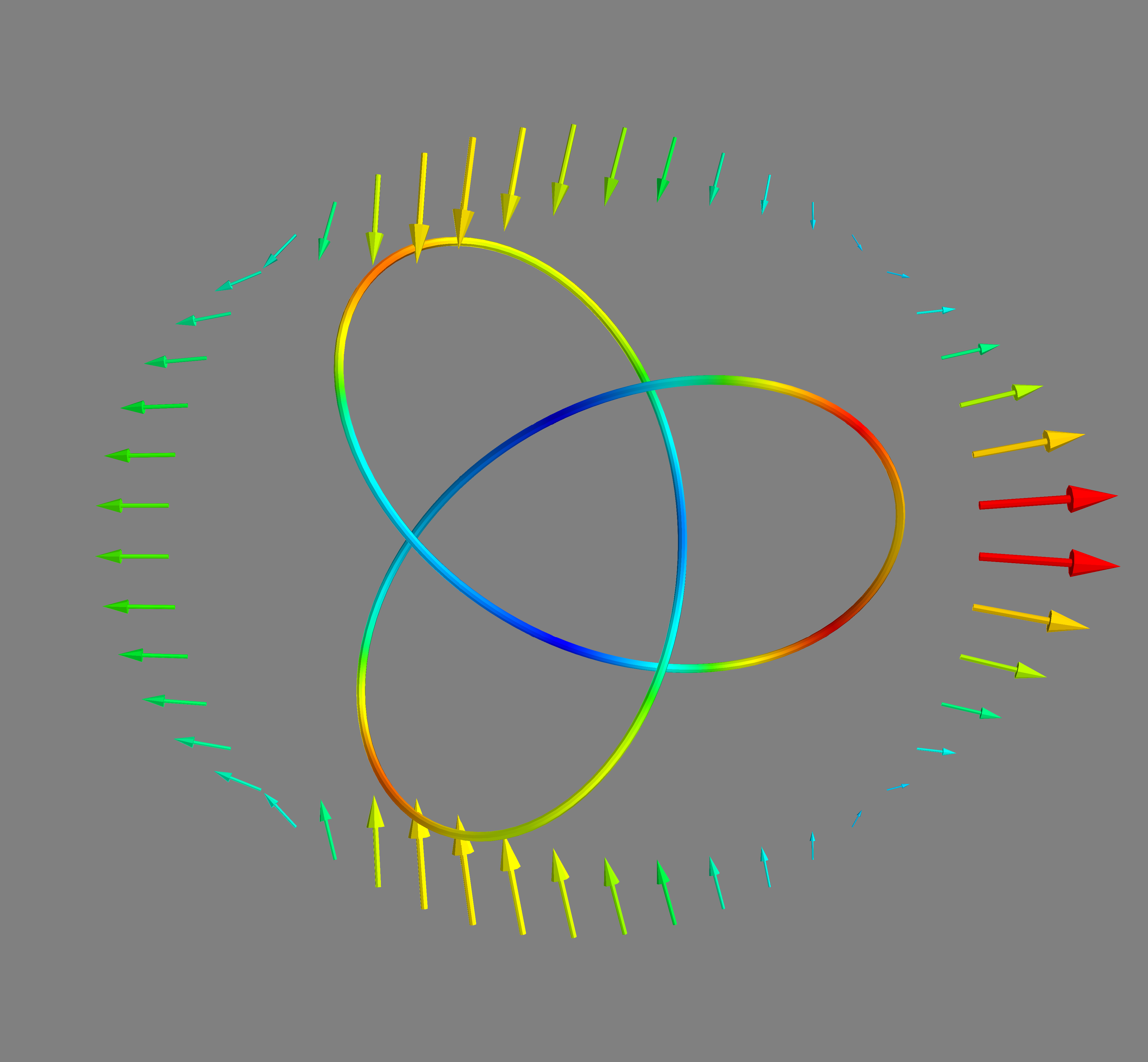}\;\;\;
\includegraphics[width=0.54\linewidth]{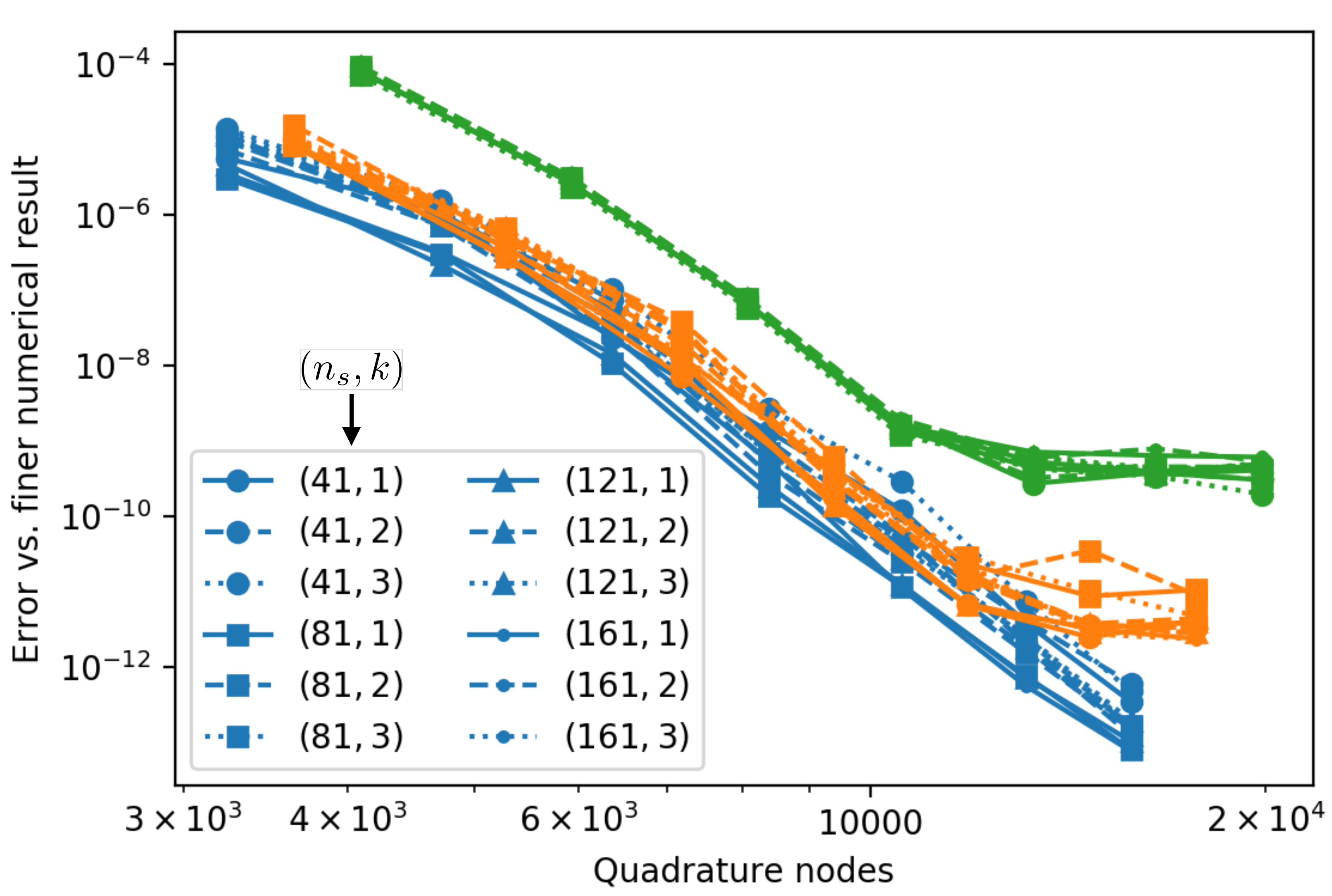}
\caption{ 
To study the convergence of our numerical method, we consider a trefoil centerline $\bm \gamma(s) =(\sin s + 2\sin 2s, \cos s - 2\cos 2s, -\sin 3s)$ and we vary the forcing frequency $k$, the fiber radius $\epsilon$, the discretization parameter $n_s$, and the number of quadrature nodes while fixing $n_\theta=7$. The flow field for $\epsilon = 0.0015$ and $k=1$ appears on the left; in this case the forcing applied is $\bm f(s) = (x(s),-y(s),0)$, and accordingly the induced fluid flow field away from the body resembles an extensional flow. The difference is that a true extensional flow would grow with distance from the center, whereas our velocity field decays with distance from the fiber. The fiber surface is colored according to the agreement of the surface velocity with Keller-Rubinow slender-body theory, with red indicating the maximum error. 
The error of our numerical method decreases exponentially in the number of quadrature nodes until it reaches a floor which depends on the fiber radius $\epsilon$. The many similar curves illustrate that the accuracy is not sensitive to the parameter $n_s$ (number of Fourier modes in the centerline direction) or the forcing wavenumber parameter $k$. The errors displayed in this figure were computed by a comparison to a reference numerical solution with quadrature parameter 39 (22496 nodes); the seven quadrature parameters displayed in the figure are evenly spaced from 15 to 33. To compare results that were computed on different periodic grids, we use Fourier interpolation to sample the centerline velocities on an evenly spaced $500$-point grid, then subtract, take the 2-norm of each of the resulting 500 vectors, and finally take the maximum ($\infty$-norm). We scale by dividing by the norm of the reference solution so that the results we report are relative. 
The blue color corresponds to $\epsilon= 0.05$, while orange and green correspond to $\epsilon = 0.005$ and $\epsilon=0.0005$. The forcing frequency parameter $k$ is equal to 1, 2, and 3, and we used $n_s = 41,81,121,161.$ Among all of these choices, the only one that seems to matter is $\epsilon$ (evidently $n_s=41$ is already sufficient to resolve the trefoil geometry). For all of these trials, the condition number of the discrete system is approximately $100/\epsilon$, that is, increasing as the fiber radius shrinks but otherwise insensitive to the discretization parameters. 
}
\label{fig:trefoilext}
\end{figure}

We now study the convergence rate of our numerical method in the presence of complicated centerline geometry and higher-frequency forcing. 
Specifically, we consider the case of a trefoil knot where the centerline and the applied force density are given respectively by 
\begin{equation}
\bm\gamma(s) = \begin{pmatrix}\sin s + 2\sin 2s\\ \cos s - 2\cos s\\ -\sin 3s\end{pmatrix}\qquad\qquad 
\bm f(s) = 
\begin{pmatrix}\sin ks + 2\sin 2ks\\ -\cos ks + 2\cos ks\\ 0\end{pmatrix}.
\end{equation}
In the case where $k=1$, the applied forcing is linear in the space variables and the resulting fluid flow resembles an extensional flow (although it decays rather than grows with distance from the fiber). This flow field is illustrated in Fig. \ref{fig:trefoilext}. We chose this trefoil curve because it has a nontrivial three-dimensional structure but does not come close to self-intersection, a more challenging problem that we will consider later in the paper.
We found that the convergence rate appears to be spectral (concave down on a log-log plot) until we reach a plateau in accuracy which depends on $\epsilon$ but not the other parameters in the problem.

\section{Comparison to other slender body theories}
\label{sec:comp2sbt}
We now compare the results of our BVP-based numerical method to the predictions of Keller-Rubinow (KR) slender-body theory as stated above \eqref{eq:sbt}. This formulation assumes that the length of the fiber has been scaled to $1$ and the parameterization has with constant speed. 
While our general numerical method does not require a unit-length fiber or even a constant-speed parameterization, in this section we numerically reparameterize all curves and then use the constant-speed, unit length parameterizations 
for both the KR slender-body theory and our method. 
The evaluation of the one-dimensional integrals in the KR version deserves some comment. The two terms in the integrand both diverge like $1/|s-t|$, but they cancel to give a bounded integrand with a jump discontinuity. In principle there is still a risk of machine arithmetic error when $s\approx t$ due to the subtraction of nearly equal quantities, but we found no problems of this kind in our experiments. For these integrals we used Gauss-Legendre quadrature with 200 nodes on each of the subintervals $[s,s+1/3]$, $[s+1/3,s+2/3]$, and $[s+2/3,s+1]$. The integrand has a jump only at $t=s$, but we also needed good resolution near $s+1/3$ and $s+2/3$ because of the particular shape of the fiber centerline in some of our tests. 

\subsection{Circular centerline, in-plane low-frequency forcing} 
\begin{figure}
\includegraphics[width=0.7\linewidth]{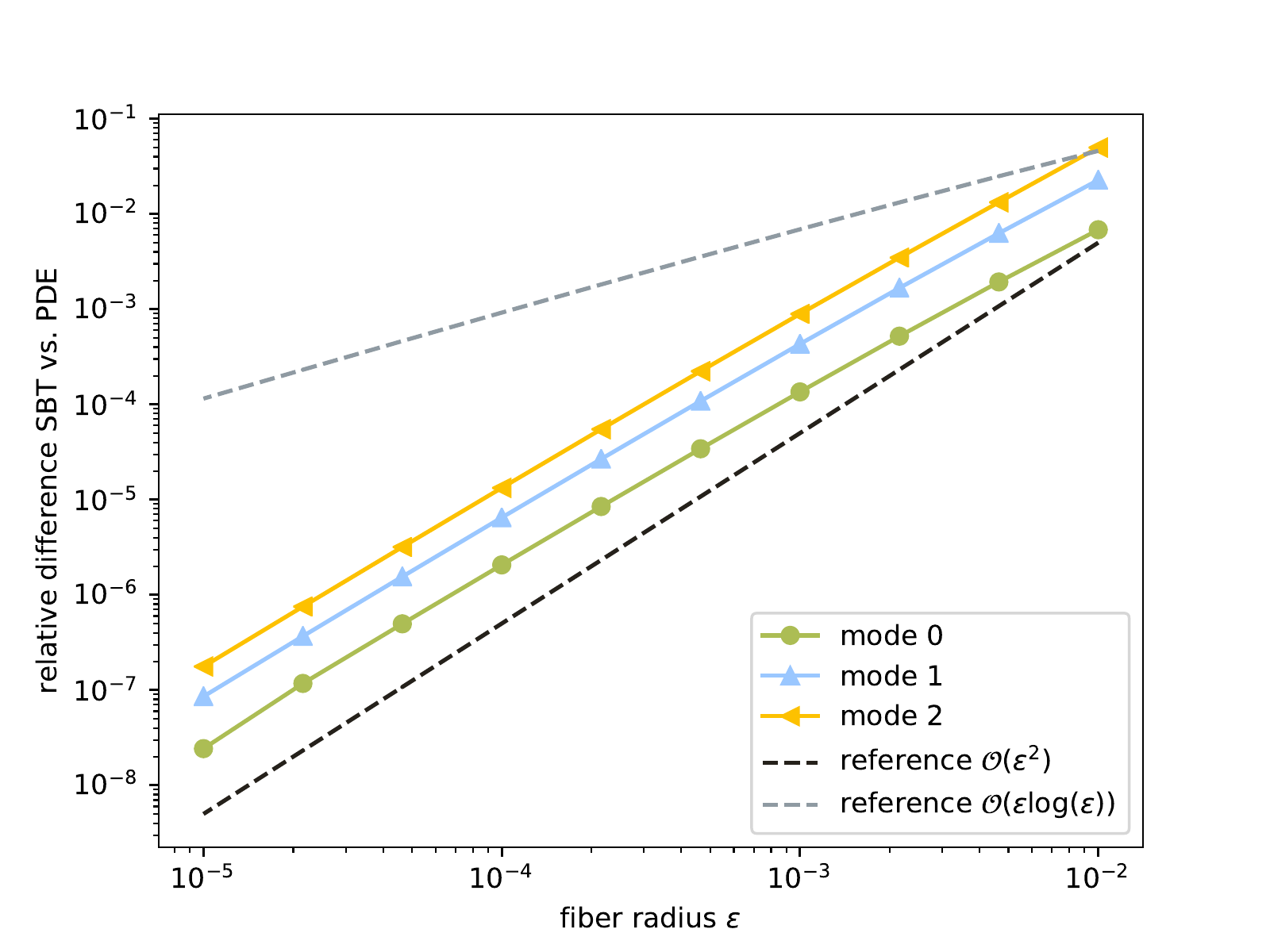}
\caption{ For a circular centerline and in-plane forcing at various wavenumbers, the discrepency between Keller-Rubinow slender body theory and our BVP-based method decays with fiber radius at rate $\mathcal{O}(\epsilon^{1.82})$. The error increases modestly with the wavenumber of the imposed centerline force density. }
\label{fig:torus_low_freq}
\end{figure}

We begin with a fiber whose centerline is circular. Instead of a uniform force density directed along the circle's symmetry axis, we consider an in-plane forcing with some sinusoidal variation: 
\begin{equation}
\bm\gamma(s) = \frac{1}{2\pi}\begin{pmatrix}\cos s \\ \sin s \\0 \end{pmatrix}\qquad\qquad 
\bm f(s) = 
\begin{pmatrix}\cos m s \\0\\0 \end{pmatrix}.
\end{equation}
The results for the first three modes $m=0,1,2$ and with $\epsilon$ ranging from $10^{-2}$ to $10^{-5}$ are given in Fig. \ref{fig:torus_low_freq}. The convergence rate is somewhat slower than $\epsilon^2$. 

\subsection{Tests with four curves}
\label{sec:fourcurves}
\begin{figure}
\includegraphics[width=\linewidth]{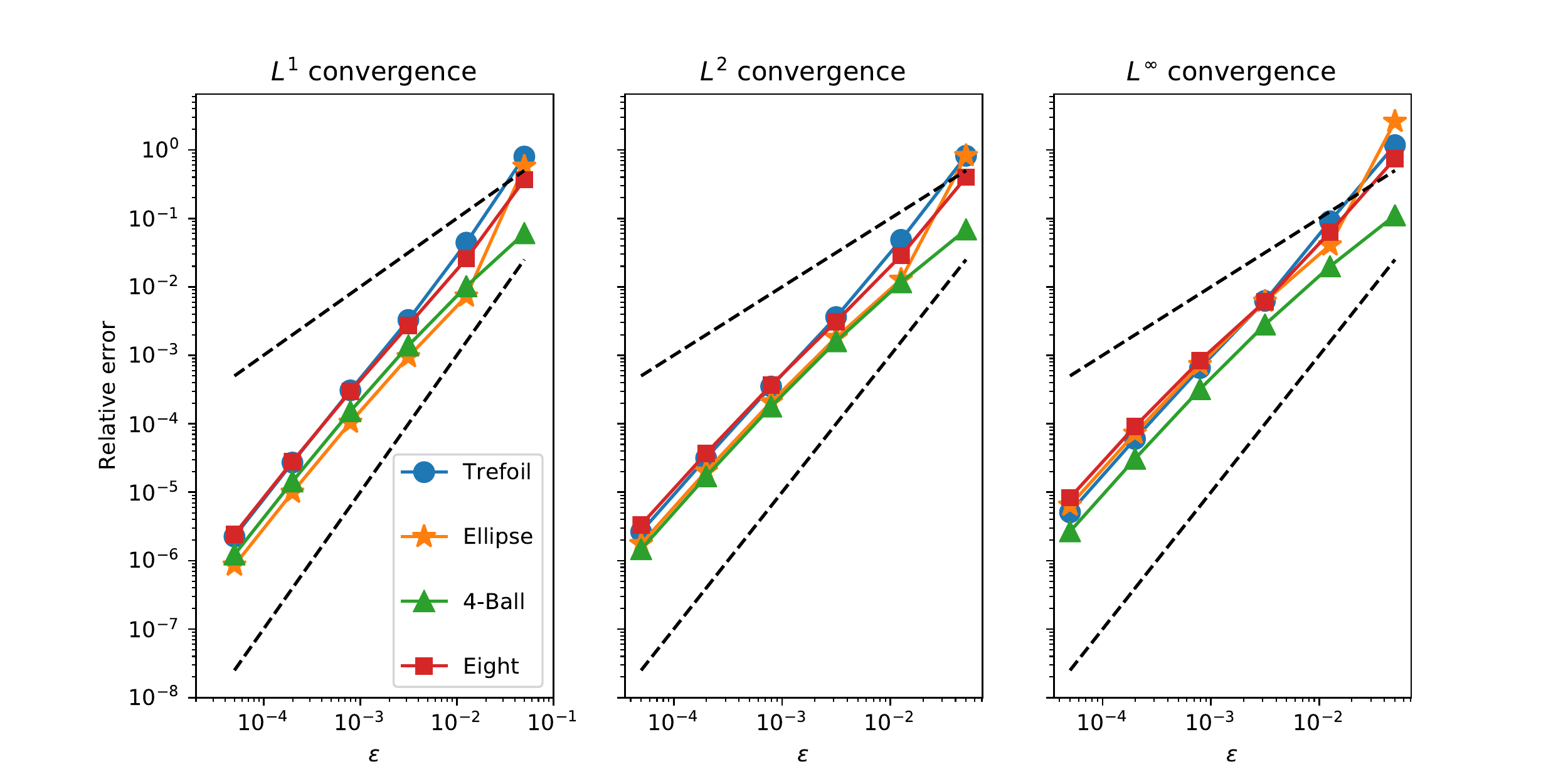}
\caption[]{The discrepancy between our BVP-based numerical method and the Keller-Rubinow slender body theory decreases at a rate between $\mathcal{O}(\epsilon)$ and $\mathcal{O}(\epsilon^2)$ (dashed reference curves) as $\epsilon\to 0$. For these calculations we used space discretization parameters $n_s = 201$ for the trefoil, ellipse, and 4-ball and $n_s=301$ for the figure-eight centerline. In all cases we took $n_\theta=5$ and quadrature parameter $n_q = 35$ (resulting in about 12,000 quadrature nodes). In reporting the relative errors we first subtract the centerline velocities of the full numerical and KR slender-body procedures, then take the 2-norm over the space dimension, then take either the $1-$, $2-$, or $\infty$-norm of the resulting vector, and finally divide by the KR slender-body norm obtained in the same way. We computed convergence rates using linear regression with the smallest three $\epsilon$-values and found that the convergence rates are comparable for all four curves and all three norms; the fastest rate was 1.78 (trefoil knot \& 1-norm) and the slowest was 1.67 (figure eight \& $\infty$-norm). 
}
\label{fig:fourcurves}
\end{figure}

We now compare Keller-Rubinow slender-body theory with our BVP-based method using several noncircular centerlines. We consider four closed curves: a planar ellipse of aspect ratio 2.5, a trefoil knot, the planar boundary of the unit ball of in the $4$-norm, that is, the curve defined by $x^4+y^4=1$, and finally a figure-eight loop. The initial parameterizations are respectively 
\begin{equation}
\bm\gamma(s) = 
\begin{pmatrix}
\cos s\\2.5\sin s\\0
\end{pmatrix},\quad
\begin{pmatrix}
\sin s + 2\sin 2s\\\cos s - 2\cos 2s\\-\sin 3s\end{pmatrix},\quad
\begin{pmatrix}
(\cos^4 s + \sin^4s)^{-1/4}\cos s\\
(\cos^4 s + \sin^4s)^{-1/4}\sin s\\0
\end{pmatrix},\quad
\begin{pmatrix}
\sin 2s\\ 1.6\sin s\\ 0.3\cos s \end{pmatrix},
\end{equation}
although we then reparameterize for constant speed and unit length as mentioned above. 
The results are shown in Fig. \ref{fig:fourcurves}. The rate at which the discrepancy decays is between $\mathcal{O}\left(\epsilon^{1.67}\right)$ and $\mathcal{O}\left(\epsilon^{1.78}\right)$ across the four centerlines and the three norms ($1$, $2$, or $\infty$).  That is, our method behaves similarly to the Keller-Rubinow formulation in these tests. 

\subsection{Exploring the limitations of SBT: self-intersection}
Keller and Rubinow's derivation is based on the method of matched asymptotic expansions; in the inner region, near the fiber surface, they assume that the fluid velocity is well approximated by the flow near a translating and rotating cylinder. This assumption is invalid when another section of the fiber is located within a distance of $\mathcal{O}(\epsilon)$. Accordingly, we expect the Keller-Rubinow to break down at least locally in the presence of near self-intersections of the fiber surface. 

To be more quantitative about what it means for the surface to nearly intersect itself, consider the quantity 
\cite{bertozzi1987extension,majda2002vorticity}
\begin{equation}\sigma(\bm\gamma) = \min_{s,t\,\in[0,1]} \frac{\|\bm\gamma(s) - \bm\gamma(t)\|}{\sin\left|\pi(s-t)\right|} .
\label{eq:sigma}
\end{equation}
Here we are assuming a constant-speed, unit length parameterization so that $\sigma$ is a geometric quantity. For a circle, $\sigma = 0.3183$, the greatest possible value. For the four curves of Section \S\ref{sec:fourcurves}, we have $\sigma = 0.0442$ for the trefoil, $\sigma=0.1738$ for the ellipse, $\sigma = 0.2826$ for the $4$-ball boundary, and $\sigma = 0.0541$ for the figure-eight loop. 

Of course, for any fixed centerline, taking a very small $\epsilon$ should result in good agreement between our method and the Keller-Rubinow formulation, just as in Fig. \ref{fig:fourcurves}. To create a test where the two methods are likely to diverge, we consider a family of centerlines where we can move the near-intersection points closer together while simultaneously reducing the fiber radius so that the ratio between the gap size and the radius remains constant. These centerlines are initially defined by
 \begin{equation}
 \bm\gamma(s) = 
 \begin{pmatrix}
  \cos(s)(1+ H  \cos3 s)\\
   \sin(s) (1+ H  \cos3 s)\\
   H\sin3 s   
\end{pmatrix} ,
\label{eq:hairtiecurve}
\end{equation}
but then numerically reparameterized and scaled so that the speed is constant and the total length is 1. 
Here $H<1$ is the distance from a point on the centerline to the unit circle (before scaling and reparameterization). When $H\to0$ the curve simplifies to a circle. When $H\to 1$, the curve points $\bm \gamma(\pi/3)$, $\bm\gamma(\pi)$ and $\bm\gamma(5\pi/3)$ all approach the origin; we refer to the distance between any two of these points, after rescaling, as the \emph{gap size}. 
Thus any value of the fiber radius $\epsilon$ exceeding half of the gap size would result in a self-intersection of the fiber surface. The relationship between $H$, the gap size, and the self-intersection quantity $\sigma$ is given in Table \ref{tbl:hairtie}. 

\begin{table}
\begin{tabular}{lcc}
\toprule 
$H$ & gap & $\sigma$\\
\midrule
0.6 &  0.05269 & 0.06084\\
0.8 &  0.02078 & 0.02400\\
0.9 &  0.00937 & 0.01082\\
0.95 & 0.00447 & 0.00515\\
0.975 &0.00218 & 0.00252\\
\bottomrule
\end{tabular}
\caption{Geometry data for the centerlines with near-self intersections defined by \eqref{eq:hairtiecurve} after reparameterization. As $H$ approaches $1$, the centerline points $\gamma(\pi/3)$ and $\gamma(\pi)$ and $\gamma(5\pi/3)$ all approach the origin. The distance between any two of these is the \emph{gap}; this approaches zero along with the self-intersection quantity $\sigma$ defined by \eqref{eq:sigma}. }
\label{tbl:hairtie}
\end{table}

\begin{figure}
\includegraphics[width = 0.6\linewidth]{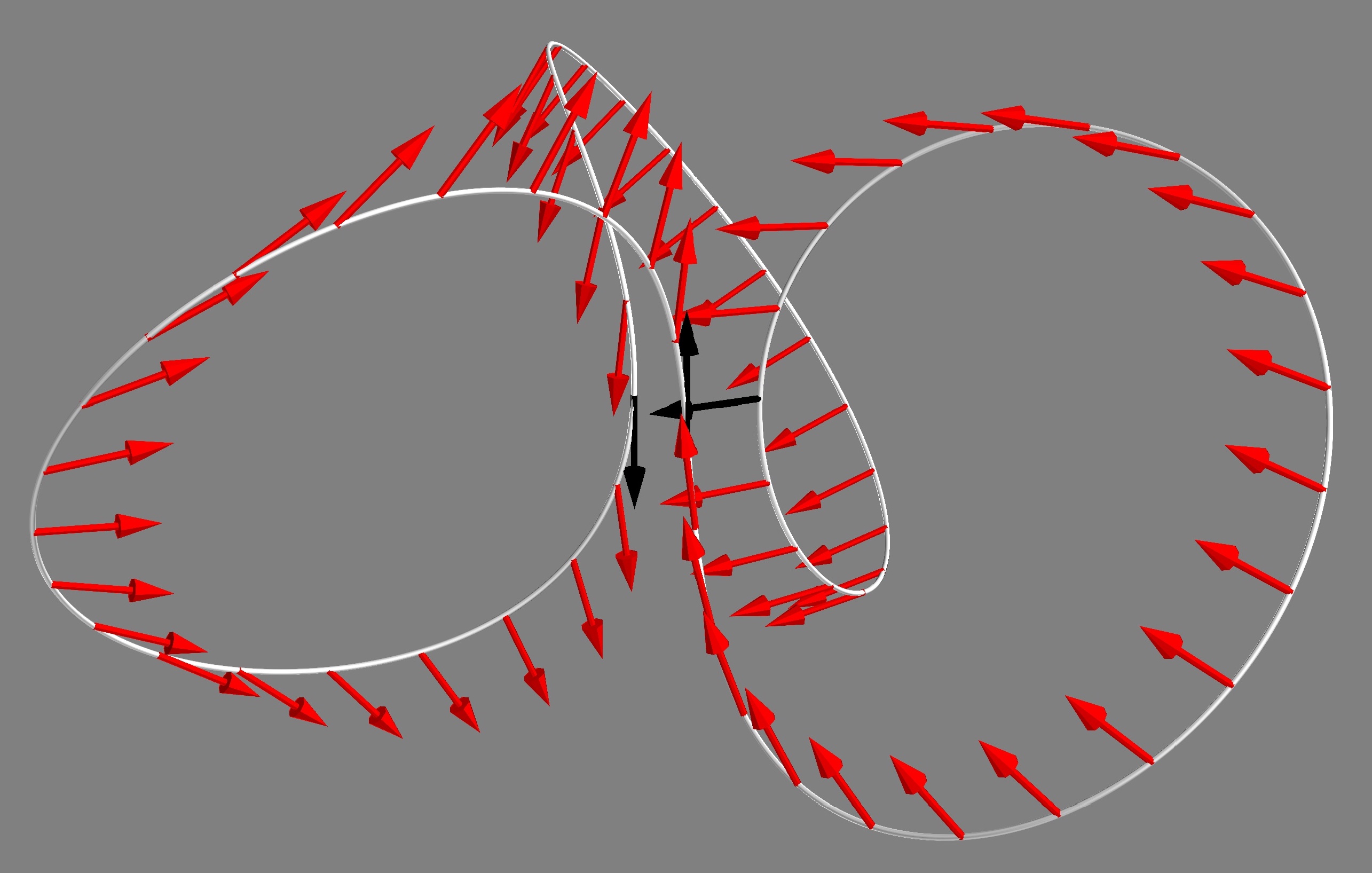}
\caption[]{The centerline defined by \eqref{eq:hairtiecurve} comes near to self-intersection when $H\approx1$; here we illustrate the curve with $H=0.8$. The arrows depict the imposed centerline force density $\bm f(s)$ given by \eqref{eq:hairtieforce}. The three black arrows indicate the forcing imposed at the centerline points $s=\pi/3$, $s=\pi$, $s=5\pi/3$ where the fiber comes closest to intersecting itself. The forcing function was chosen so that these three branches of the fiber are pushed in contrasting directions, specifically $\bm{\hat{z}}$, $\bm{\hat{x}}$, and $-\bm{\hat{z}}$. }
\label{fig:sbtfail_arrows}
\end{figure}

As an example, Fig. \ref{fig:sbtfail_arrows} shows the centerline that results from setting $H=0.8$. The arrows indicate the centerline force function we impose, given by
\begin{equation}\bm f(s) = \begin{pmatrix}
-\cos\left(s + \frac{\pi}{3\sqrt{3}}\sin(s)\right)\\
0\\
\sin\left(s + \frac{\pi}{3\sqrt{3}}\sin(s)\right)
\end{pmatrix} .
\label{eq:hairtieforce}
\end{equation}
We chose this peculiar form for the forcing function because of its values at the three points of nearest self-intersection, illustrated as black arrows in Fig. \ref{fig:sbtfail_arrows}: $\bm f(\pi/3) = \bm{\hat{z}}$, $\bm f(5\pi/3) = -\bm{\hat{z}}$, and $\bm f(\pi) = \bm{\hat{x}}$. That is, two of the three branches that pass near the origin are being forced in opposite directions tangential to the fiber surface while the third branch is being forced in the plane normal to the centerline. This complicated scenario is designed to find breakdowns in the Keller-Rubinow slender-body theory. 

\begin{figure}
\includegraphics[width = 0.9\linewidth]{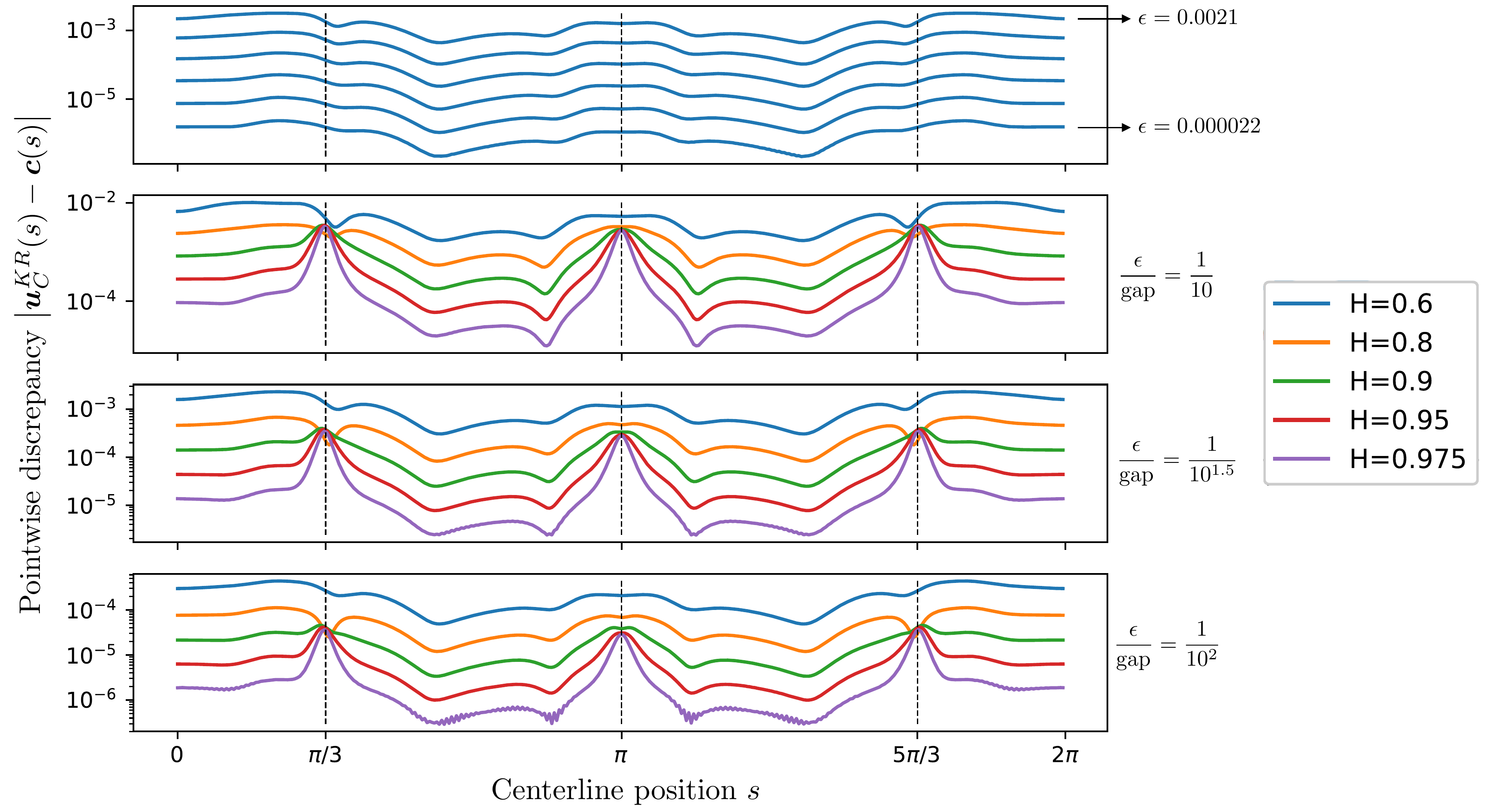}
\caption[]{The breakdown of Keller-Rubinow slender body theory when the fiber surface approaches itself is local. In the top panel, the centerline is fixed while the fiber radius decreases; the discrepancy between the KR centerline velocity $\bm u_C^{KR}(s)$ and our centerline velocity function $\bm c(s)$ decreases uniformly in $s$.  In the lower three panels, the centerline shifts ($H\to 1$ in \eqref{eq:hairtiecurve}) and the radius decreases simultaneously so that the ratio of the radius to the gap size remains constant at $1/10$, $1/10^{1.5}$, or $1/10^2$. In these three cases we see that the pointwise discrepancy decreases along most of the fiber length as $\epsilon\to0$. However, the errors are $\mathcal{O}(1)$ at the three locations $s=\pi/3$, $s=\pi$, and $s=5\pi/3$ where the fiber comes close to intersecting itself. }
\label{fig:sbtfail_pw}
\end{figure}

\begin{figure}
\includegraphics[width = 0.6\linewidth]{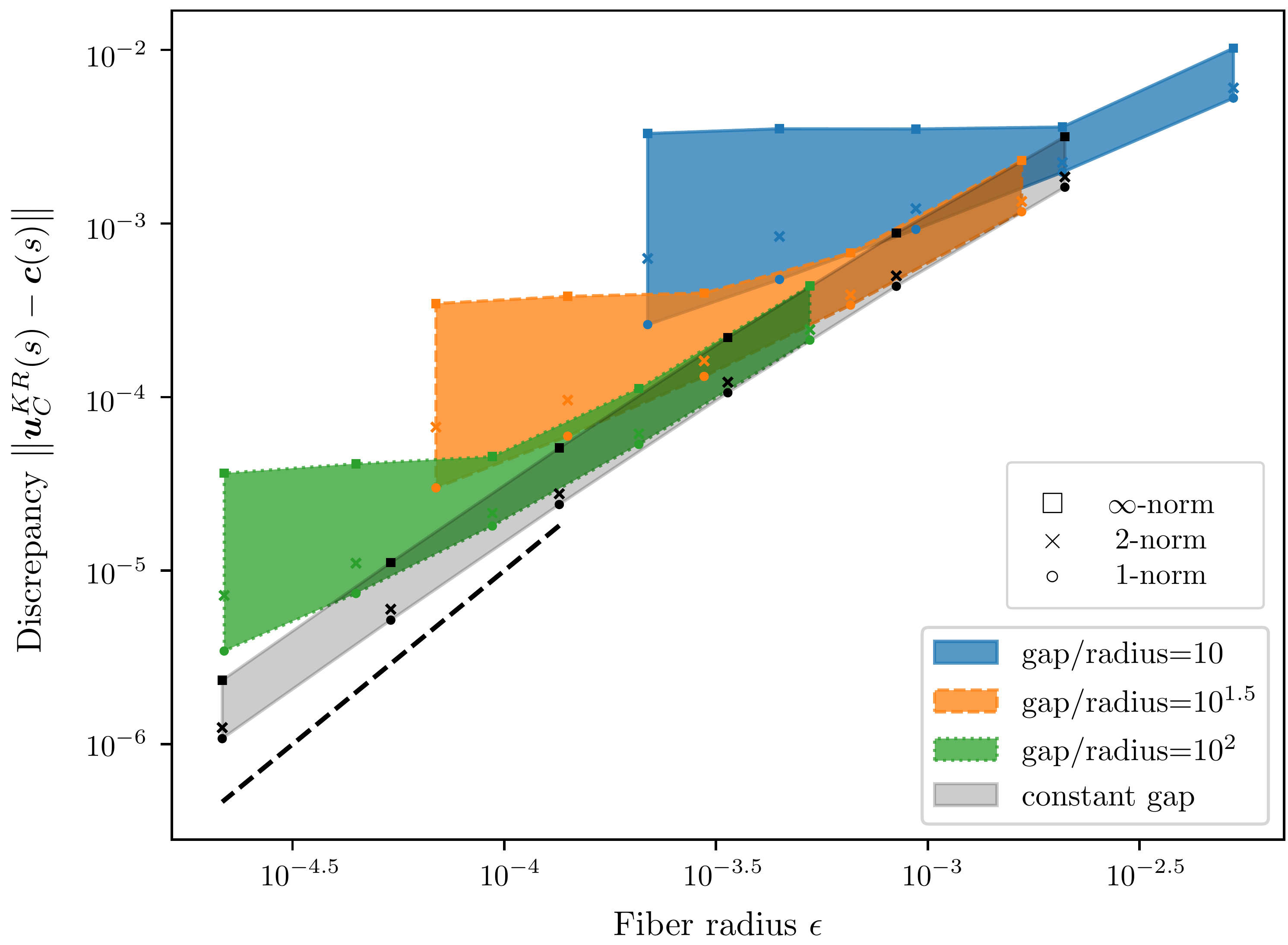}
\caption[]{The gray polygon shows the convergence of the Keller-Rubinow slender body theory and our BVP-based method when the centerline is fixed and the fiber radius shrinks, at about the same order we observed in the other tests, $\mathcal{O}(\epsilon^{1.69})$. The dashed reference curve is $\mathcal{O}(\epsilon^2)$. The upper boundary of the gray polygon gives the $\infty$-norm while the lower boundary is the $1$-norm; the $2$-norm values appear as $\times$ markers inside the polygon. The similarity of the convergence in the $1$- and $\infty$-norms suggests that the convergence must be uniform in the centerline coordinate $s$, as is the case (top panel of Fig. \ref{fig:sbtfail_pw}). The blue, yellow, and green colored polygons describe tests wherein the fiber centerline approaches itself while its radius decreases simultaneously. In these cases the $\infty$-norm stagnates while the $1$-norm continues to decrease, suggesting that the breakdown of the Keller-Rubinow formulation is local in $s$ (lower three panels of Fig. \ref{fig:sbtfail_pw}). }
\label{fig:sbtfail_hockey}
\end{figure}

Figures \ref{fig:sbtfail_pw} and \ref{fig:sbtfail_hockey} illustrate the breakdown of the Keller-Rubinow formulation with self-intersections. In one set of simulations we fix $H=0.6$ and we let $\epsilon\to0$. As expected, our method agrees with the Keller-Rubinow formulation at approximately order $\mathcal{O}(\epsilon^{1.7})$ (gray polygon in Fig. \ref{fig:sbtfail_hockey}) and this convergence is uniform in the centerline coordinate $s$ (top panel of Fig. \ref{fig:sbtfail_pw}). Indeed, the fact that the convergence rate is essentially the same in the $1$-, $2$- and $\infty$-norms suggests that the discrepancy between the two methods must decay uniformly in $s$. This situation changes, however, when we allow the parameter $H$ to increase toward $1$ while keeping the fiber radius at a constant fraction of the gap size. We did this with ratios of $1/10$ (second row in Fig. \ref{fig:sbtfail_pw}, blue polygon in Fig. \ref{fig:sbtfail_hockey}), $1/10^{1.5}$ (third row in Fig. \ref{fig:sbtfail_pw}, orange polygon in Fig. \ref{fig:sbtfail_hockey}), and $1/100$ (fourth row in Fig. \ref{fig:sbtfail_pw}, green polygon in Fig. \ref{fig:sbtfail_hockey}). In all of these cases we see that the local error does not decrease uniformly. Near the intersection points at $s\in\{\pi/3,\pi,5\pi/3\}$ the two methods appear to be converging to different answers and the discrepancy is $\mathcal{O}(1)$. At fiber points away from the intersection regions, the discrepancy does decay at about the same rate. Thus, the overall error stagnates when measured in the $\infty$-norm (the upper boundary of the polygons in Fig. \ref{fig:sbtfail_hockey}) but continues to decrease in the $1$-norm (the lower boundary of the polygons in Fig. \ref{fig:sbtfail_hockey}). The breakdown in the Keller-Rubinow formulation follows from the invalidity of their inner expansion of velocity when another part of the fiber surface lies within a distance of $\mathcal{O}(\epsilon)$. It is interesting that we can detect this breakdown even when the other fiber surface is relatively far away from the observation point ($100\epsilon$ in the last set of tests).

\section{Conclusion and future work}
In this paper we presented evidence that the standard slender body theories based on matched asymptotic expansions can give inaccurate results when the fiber surface approaches itself, a situation that is common in biophysical and industrial processes. The breakdown is local in the sense that these formulations still give good results on isolated sections of the fiber. As an alternative, we gave a numerical method for a fully three-dimensional slender body Stokes boundary value problem which can be stated with reference to a one-dimensional centerline force density only. 

Some immediate next steps would be to account for multiple fibers \cite{li2013sedimentation,hamalainen2011papermaking}, fibers with free ends \cite{free_ends,johnson1980improved}, and dynamic problems \cite{shelley2000stokesian,lindner2015elastic}. Another possible extension of this work would be to consider inertial flows, where the underlying PDE is different but the fiber integrity and average force conditions at the boundary are still a reasonable way to make physical sense of one-dimensional forcing data. 

While our method is more accurate in the presence of near contacts, it is much more computationally intensive than the Keller-Rubinow formulation. For example, the largest simulations reported here required assembling and solving a dense system of 8712 linear equations, wherein the majority of the matrix entries are difficult integrals requiring two-dimensional quadrature. To simulate many fibers or to move from static to dynamic problems, we will need some combination of a cheaper algorithm, larger machines and/or parallelism. At larger scales, we will also need to reconsider our choice to use the full SVD for the numerical linear algebra. A fast multipole method may be an appropriate strategy for larger sized problems \cite{ROKHLIN1985187}. 

Given the adequate performance of the inexpensive Keller-Rubinow formulation for isolated sections of the fiber, a hybridization of the two methods might provide another way to reduce the computational expense relative to the method presented here. 
\label{sec:conclusion}

\section{Acknowledgements}
WM and HB were supported by NSF grant DMS-1907796. DS was supported in part by NSF grant DMS-2009352. LO was supported by NSF Postdoctoral Research Fellowship DMS-2001959. YM was supported by NSF (DMS-2042144) and the Math+X grant from the Simons Foundation.
We acknowledge the generous hospitality and computational resources of the IMA, where the project was initiated. 
We thank the anonymous reviewers for their careful reading of the manuscript and for their helpful insights.

\bibliographystyle{plain}
\bibliography{NotesNotes}

\appendix
\section{The Bishop frame and fiber surface} 
Here we give more specific information on our surface geometry construction using a quaternion-based initial value problem. We first describe the continuous formulation and then give numerical details. 
\subsection{Defining the Bishop frame using an initial value problem}
Let $\bm \gamma: [0,2\pi]\to \mathbb{R}^3$ be a twice-differentiable and periodic curve parameterizing the centerline of a closed fiber. We do not require an arclength parameterization, that is, we require only $|{\bm\gamma}'|>0$ instead of $|{\bm\gamma}'|=1$.  The unit tangent vector and its derivative with respect to $s$ are defined by 
\begin{align}
\label{eq:unitT}\bm T(s) &= \frac{1}{|{\bm \gamma}'|}{\bm \gamma'} , \\
{\bm T'}(s) &= \frac{1}{|{\bm \gamma}'|}{\bm \gamma}'' - \frac{{\bm \gamma}'\cdot {\bm \gamma}''}{|{\bm \gamma}'|^3}{\bm \gamma}' .
\end{align}
Because the vector $\bm T$ is perpendicular to its derivative, it traces out a closed loop on the unit sphere as $s$ ranges from $0$ to $2\pi$. We wish to find vectors $\bm N(s)$ and $\bm B(s)$ which together with $\bm T(s)$ form an orthonormal frame.  The simplest method, due to Frenet, is to set $\bm N = {\bm T}' / |{\bm T}'|$ and $\bm B = \bm T\times\bm N$; these definitions are \emph{local} in the sense that we need to know $\bm\gamma$ and its derivatives only at a fixed $s$ in order to determine the frame $\{\bm T(s),\bm N(s),\bm B(s)\}$. However, the Frenet normal and binormal vectors are undefined wherever ${\bm T}' = \bm 0$, or equivalently whenever the acceleration $ {\bm \gamma}''$ is a scalar multiple of the velocity ${\bm \gamma}'$. A far-from-pathological example where this occurs is the boundary of the unit ball in the $\ell^4$ norm in the plane, which has the polar parameterization $r(s) = \left(\cos^4s+\sin^4s\right)^{1/4}$. 
We therefore use an alternative \emph{global} frame which is defined as the solution to a certain initial value problem. 
Let $\bm {T_0} = \bm T(0)$ and let $\bm N_0$ and $\bm B_0$ be any vectors completing the frame at $s=0$.  
Then define $\bm\omega = \bm\omega(s)$ by 
\begin{equation}
\label{eq:bishopomega}
\bm\omega = \bm T\times {\bm T}' + \frac{|{\bm \gamma}'|\alpha}{L}\bm T
\end{equation}
where $\alpha$ is a constant scalar to be determined later and $L$ is the total length of the path $\bm \gamma$. 
Now consider the differential equation of rotation
\begin{equation}
\label{eq:bishopdiffeq}
{\bm v}' = \bm\omega \times \bm v.
\end{equation}
If the initial condition is $\bm v_0 = \bm T_0$, we find that the solution is $\bm v(s) = \bm T(s)$; to see this, use the triple cross product formula $(\bm A\times\bm B)\times \bm C = (\bm C\cdot\bm A)\bm B - (\bm C\cdot\bm B)\bm A$ to write
\begin{equation}\left(\bm T\times {\bm T}' + \frac{\alpha}{2\pi}\bm T\right)\times \bm T = (\bm T\cdot\bm T){\bm T}' -(\bm T\cdot{\bm T}'){\bm T} + 0 = {\bm T}' .\end{equation}
This is a global definition for $\bm T$ which is identical to the local version \eqref{eq:unitT}. However, 
we can also solve \eqref{eq:bishopdiffeq} with the initial conditions $\bm v_0 = \bm N_0$ and $\bm v_0 = \bm B_0$ to obtain unit vector functions $\bm N(s)$ and $\bm B(s)$. These complete the frame for each $s\in[0,2\pi],$ for if $\bm v_1$ and $\bm v_2$ are solutions of \eqref{eq:bishopdiffeq}, we have $\frac{d}{ds}(\bm v_1 \cdot \bm v_2) = \bm v_1\cdot(\bm \omega\times \bm v_2) + (\bm\omega\times \bm v_1)\cdot \bm v_2 = 0$, showing that evolution under \eqref{eq:bishopdiffeq} does not change lengths or angles between vectors. 

The parameter $\alpha$ determines the speed of rotation about $\bm T$ and therefore affects the evolution of $\bm N$ and $\bm B$ but not $\bm T$. We require a value of $\alpha$ giving a periodic frame: $\bm N(2\pi) = \bm N(0)$ and $\bm B(2\pi) = \bm B(0)$.  This value of $\alpha$ is not unique; one can always add an integer to obtain another periodic frame where $\bm N$ and $\bm B$ twist around $\bm T$ a different number of times on $[0,2\pi]$.  One method of choosing $\alpha$ is to solve \eqref{eq:bishopdiffeq} with $\alpha =0$ and then examine the rotation carrying $\{\bm N(2\pi), \bm B(2\pi)\}$ to $\{\bm N_0, \bm B_0\}$; then we set $\alpha$ to be the angle by which the initial and final $\bm N$ and $\bm B$ differ and solve \eqref{eq:bishopdiffeq} again.  

By defining the scalar functions $\kappa_1 = \frac{1}{|{\bm \gamma}'|} {\bm T}'\cdot \bm N$, $\kappa_2 = \frac{1}{|{\bm \gamma}'|} {\bm T}'\cdot \bm B$, $\kappa_3 = \frac{\alpha }{L}$, we can make the differential equation \eqref{eq:bishopdiffeq} equivalent to the usual formulation of the Bishop frame, which is to start with arbitrary $\kappa_i$ and evolve the differential equation
\begin{equation}
\label{eq:tradbish}
\frac{d}{ds}\begin{pmatrix}\bm T\\\bm N\\\bm B\end{pmatrix}
= 
|\bm\gamma'|\begin{pmatrix}0&\kappa_1&\kappa_2\\-\kappa_1&0&\kappa_3\\
-\kappa_2&-\kappa_3&0
\end{pmatrix}
\begin{pmatrix}\bm T\\\bm N\\\bm B\end{pmatrix}.
\end{equation}
There is a factor of $|{\bm \gamma}'|$ on the right side of \eqref{eq:tradbish} because we have not assumed an arclength parameterization. 

We can now use the centerline curve $\bm\gamma$ and the orthonormal frame $\{\bm T,\bm N,\bm B\}$ to define the fiber surface: let  
\begin{equation}
\bm X(s,\theta) = \bm \gamma(s) + \epsilon\cos(\theta)\bm N(s) + \epsilon \sin(\theta)\bm B(s).
\label{eq:fibersurfaceNB}
\end{equation}
The Jacobian of this transformation and the normal vector pointing out of the fiber into the fluid can be calculated directly from $\bm X_\theta \times \bm X_s$ by writing all vectors in the right-handed basis $\{\bm T,\bm N,\bm B\}$. The surface normal is 
\begin{equation}
\bm \nu(s,\theta) = \cos(\theta)\bm N(s) + \sin(\theta)\bm B(s)
\end{equation}
and the area element is 
\begin{equation}
J(s,\theta) = \epsilon |{\bm \gamma}'|\Big(1-\epsilon \cos(\theta)\kappa_1(s) - \epsilon\sin(\theta)\kappa_2(s)\Big).
\end{equation}

\subsection{Quaternion-based numerical calculation of the Bishop frame}
We will focus on the case where the centerline $\bm \gamma$ is defined as the Fourier interpolant of some discrete set of points in $\mathbb{R}^3$ rather than by a symbolic formula. In this setting we wish to compute the frame vectors and their derivatives accurately. To reduce the number of independent quantities which must be numerically integrated from six (for $\bm N$ and $\bm B$, since $\bm T$ is known independently) to four, we use a formulation in terms of quaternions. We write quaternions as vector-scalar pairs: $q = (\bm z,r)$ is the same as $q = z_1 i +z_2j + z_3 k + r$. 
Then we employ a reformulation of \eqref{eq:bishopdiffeq} derived by \cite{boyle2017integration}:
\begin{equation}
 q' = \frac{1}{2}(\bm\omega,0) * q .
\label{eq:bishopquaternion}
\end{equation} 
Here the symbol $*$ denotes the noncommutative multiplication of quaternions, and a reformulation of \eqref{eq:bishopquaternion} using the vector dot and cross products is 
\begin{equation}
 q' = \frac{1}{2}\left( {\bm z}',r'\right) =  \left(r\bm\omega + \bm\omega\times \bm z, -\bm \omega\cdot\bm z  \right)  .
\label{eq:bishopcross}
\end{equation}
We then solve \eqref{eq:bishopquaternion} with the initial condition $q_0 = ((0,0,0),1)$ to obtain a solution as a path in $\mathbb{R}^4$. From this we can obtain $\bm T$ and $\bm N$ and $\bm B$ from 
\begin{align}
\label{eq:globalTquat}
(\bm T(t),0) &= q(t) * (\bm T_0,0) * q(t)^{-1}\\
(\bm N(t),0) &= q(t) * (\bm N_0,0) * q(t)^{-1}\\
(\bm B(t),0) &= q(t) * (\bm B_0,0) * q(t)^{-1}.
\end{align}
This formulation has the advantage that the computed frame $\{\bm T, \bm N, \bm B\}$ is exactly a rotated version of $\{\bm T_0, \bm N_0, \bm B_0\}$, even in the presence of numerical errors in the computation of $q(t)$.  The disadvantage is that numerical error in $q(t)$ may cause the globally computed $\bm T$ to differ from the local one, i.e. \eqref{eq:globalTquat} may not be identical to \eqref{eq:unitT}. 
One may assess the accuracy of the integration by comparing the two versions. 
Furthermore, one can take advantage of the local method of computing $\bm T$ by adding an auxiliary rotation spinning the global version of $\bm T$ towards the local one.  This amounts to adding a term to \eqref{eq:bishopomega}: 
\begin{equation}
\label{eq:bishopomega2}
\bm\omega(t) = \bm T(t)\times {\bm T}'(t) + \frac{\alpha}{2\pi}\bm T(t) + 10 \bm T^{num}(t) \times \bm T(t).
\end{equation}
In \eqref{eq:bishopomega2} the vector $\bm T$ is computed from the local definition \eqref{eq:unitT} while the vector $\bm T^{num}$ is the global version \eqref{eq:globalTquat}. Typically these are nearly identical and so their cross product is small, i.e. \eqref{eq:bishopomega} and \eqref{eq:bishopomega2} are very similar. The constant factor of $10$ is a heuristic which improves the accuracy of the computation without significantly increasing the stiffness of the problem. 
We use the SciPy integration routine {\tt solve\_ivp} \cite{2020SciPy-NMeth} to solve \eqref{eq:bishopquaternion} with both absolute and relative error tolerances set to $5\cdot10^{-14}$.

\section{Quadrature for single- and double-layer potentials on a slender fiber}
\label{sec:quadratureappendix}
\begin{figure}
\begin{tabular}{m{0.64\linewidth}m{0.32\linewidth}}
\includegraphics[width=\linewidth]{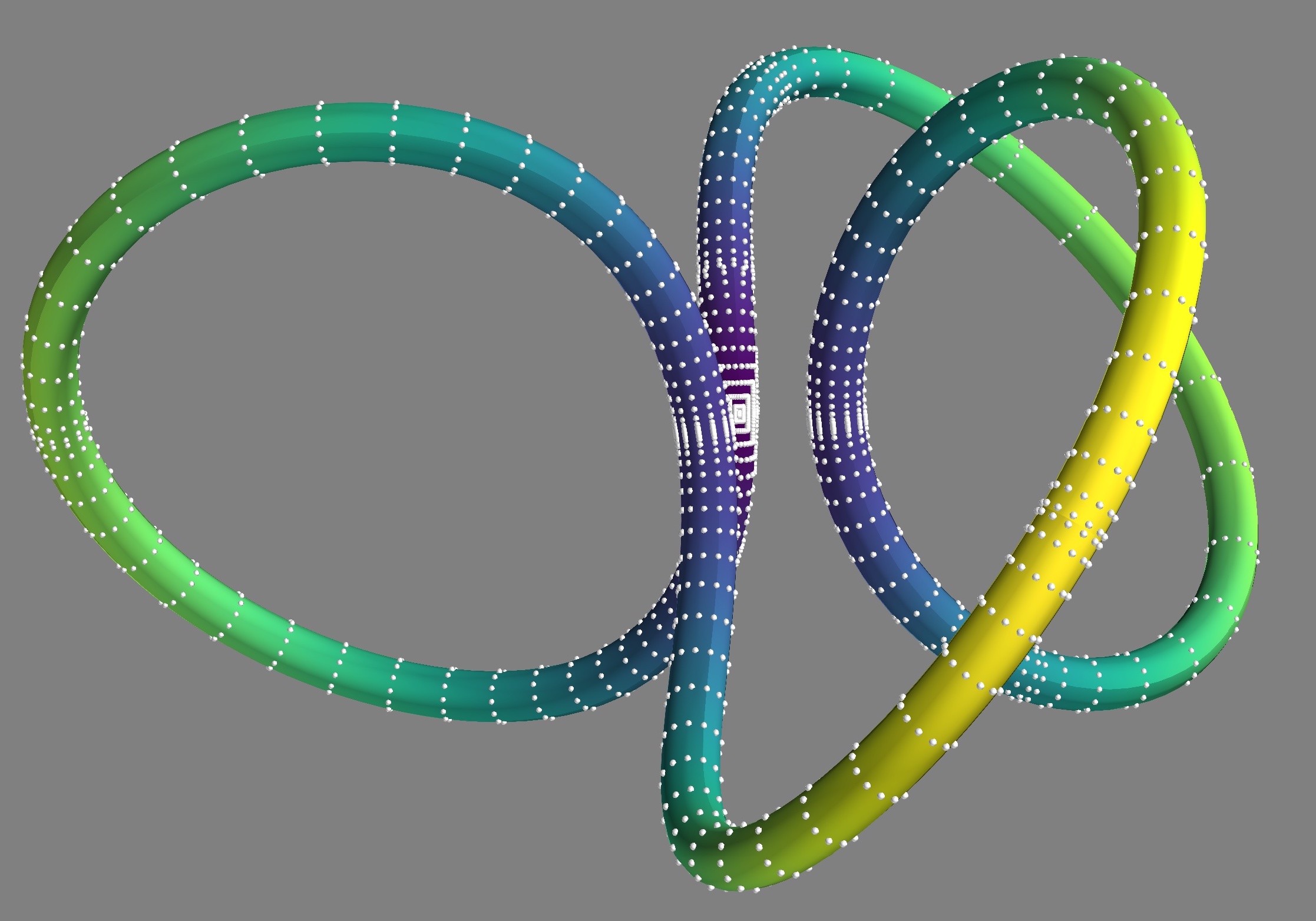}\,
&
\includegraphics[width=\linewidth]{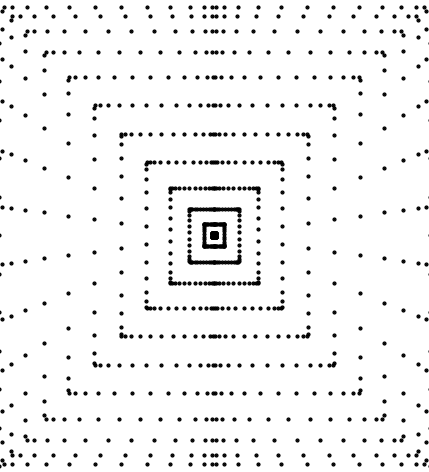}
\end{tabular}
\vspace{0.2in}
\includegraphics[width=\linewidth]{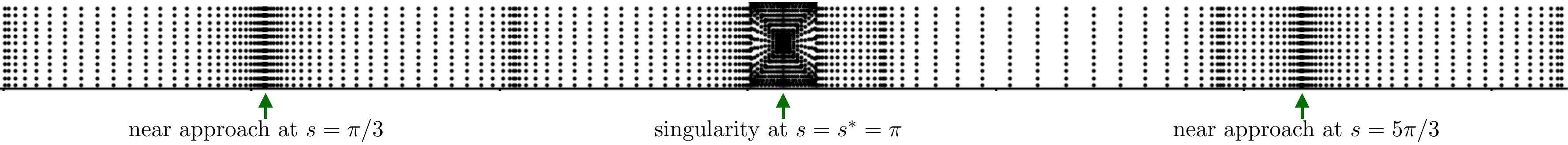}\\
\caption{The quadrature is designed to integrate a singularity of order $1/r$ at a source point. Above, we show the rule generated for a fiber surface which approaches the origin three times, specifically for the centerline \eqref{eq:hairtiecurve} with $H=0.8$ and $\epsilon$ equal to one fifth of the gap size. The color on the fiber surface indicates distance from the source point, or more precisely $|\bm \gamma(s)-\bm\gamma(s^*)|$. The source point $(s^*=\pi, \theta^*=0)$ is visible at the center of the figure on the rear fiber branch. The quadrature nodes cluster in nested squares around the singularity, as the right-hand figure indicates. Below this we plot the same quadrature nodes in the flat $s\theta$-plane. The innermost region is decomposed into six triangles and the outer region is broken into rectangles where we use rescaled Gauss-Legendre-trapezoid grids. If the singularity moves in the circumferential $(\theta)$ direction, the quadrature rule can be periodically shifted and reused; for different centerline values $s^*$ the rule must be regenerated. }
\label{fig:sixtri}
\end{figure}

We now construct quadrature rules for the single- and double-layer potentials on thin tubes. Let the centerline $\bm \gamma(s)$ and fiber surface $\bm X(s,\theta)$ be constructed as above. If the integrand were smooth, the doubly periodic geometry would be well suited to quadrature on a regular grid (double trapezoid rule). However, we have not found a way to take advantage of the double periodicity for the singular integrands required in the velocity and traction integrals  \eqref{eq:dv}-\eqref{eq:df}. Instead, we decompose the fiber surface into subregions by $s$. The inner region consists of all $s$ satisfying $|\bm\gamma(s)-\bm\gamma(s^*)| < 5\epsilon$; this is then further split into six triangles, each of which has a vertex at the source point. 
The outer region is divided into subrectangles where we use an exponentially rescaled Gauss-Legendre integration in the $s$-direction and the trapezoid rule in the $\theta$-direction, thus taking advantage of periodicity in one of the two directions. The procedure accounts for situations where the fiber centerline approaches itself and resolves the integrand more finely as needed. An example of the resulting quadrature nodes is given in Fig. \ref{fig:sixtri}. Our quadrature generation procedure relies on knowledge of the fiber centerline geometry and the critical centerline value $s^*$; thus the rules must be regenerated for each new centerline and even for each new value of $s^*$. However, they can be shifted periodically and reused when $s^*$ is fixed but $\theta^*$ changes. The accuracy and size of the quadrature is controlled by a single parameter $q_n$ and the total number of nodes is $\mathcal{O}(q_n^3)$. The results presented in this paper used values of $q_n$ ranging from $30$ to $45$. We give a more detailed description in the following subsections. 

\subsection{Inner region: integration of the singularity using triangles} 
We begin by describing a quadrature rule for the square $[-1,1]^2$ which is well adapted for a $1/r$ singularity with angular dependence at the origin.   
Consider first the triangle $T$ with vertices at $(0,0)$, $(0,1)$ and $(1,1)$. An integral over $T$ can be transformed to an integral on the square $[0,1]^2$ via 
\begin{equation}
\int_T\psi\,dA = \int_{0}^1\int_{0}^1 \psi\left(uv,u\right)u \,dv\,du .
\end{equation}
This mapping expands the region around the singularity, providing a regularizing factor which yields a smooth and bounded integrand on the square. We use an outer product of Gauss-Legendre integration rules in the $uv$-domain, which leads to clustering of the nodes around the singularity in the triangular domain. 

We use similar strategies for the other five subtriangles depicted on the right panel of Fig. \ref{fig:sixtri}. This leads to a quadrature on the square $[-1,1]^1$.  
Finally, we shift and scale the square to the inner (singular) integration region in the $s\theta$-plane, that is, with $\theta\in(\theta^*-\pi , \theta^*+\pi)$ and $s$ satisfying $|\bm\gamma(s)-\bm\gamma(s^*)| < 5\epsilon$. The right and left halves are scaled differently in $s$ if $|\bm\gamma'|$ is not constant near $s^*$. 

The total number of quadrature nodes in this inner region is $6q_n^2$. We have not carefully optimized this quadrature for the inner region because the outer region generally requires more nodes.

\subsection{Outer region: integration over multiple scales using rectangles} 

It remains to integrate over a region where the integrand is smooth, that is, the part of the fiber surface where $|\bm\gamma(s)-\bm\gamma(s^*)| > 5\epsilon$. We use an outer product of one-dimensional rules: in the circumferential coordinate $\theta$ we use the trapezoid rule, while for the centerline coordinate $s$ we first subdivide and then use an exponentially rescaled Gauss-Legendre integration on each subinterval.  

The subdivision in $s$ allows us to respond to the fact that the integration on the outer region becomes more challenging as $\epsilon\to0$. Indeed, the length of the outer domain in $s$ is approximately $2\pi - 10\epsilon|\bm \gamma'(s^*)|$, and the singularity therefore lies at a distance of $5\epsilon|\bm \gamma'(s^*)|$ from either end of the integration interval; thus as $\epsilon\to0$ the integration region approaches the singularity more closely. Therefore, when $\epsilon$ is small we want to put more subintervals near the ends of the integration region. 

The possibility of near self-intersection of the fiber centerline presents another integration challenge. In this situation we may have a large increase in the value of the integrand inside a very narrow subregion which could appear anywhere on the subinterval. We would like to concentrate some small integration subpanels in such self-intersection zones. 

We propose a method which addresses these issues, at the cost of becoming somewhat complicated and specializing to a specific fiber centerline and a specific value of $s^*$. We begin by placing subdivision markers at any value of $s$ where $h(s) = |\bm \gamma(s) - \bm\gamma(s^*)|$ has a local extremum. Additionally, we place subdivision markers at any location where $\log_{10}(h(s) / (5\epsilon))$ has a positive integer value. This results in a concentration of fiber subintervals in any location where the centerline approaches $\bm\gamma(s^*)$, including the ends of the outer integration region; the concentration becomes more extreme as the fiber radius shrinks. However, if there are many local extrema this procedure may produce too many subdivisions in places where the integrand is relatively smooth. Thus, we remove the subdivision markers at the extrema unless a centered-difference estimation of $\frac{d^2}{ds^2}\log_{10}(h(s))$ gives a relatively large value (greater than four times the average on a regular 30-point grid).  

Once the subdivisions are finalized, we create outer-product rules for each subrectangle. We always use $n_q$ trapezoidal nodes in the periodic $\theta$-direction, while the number of nodes in the $s$-direction varies with the length of the subinterval and the logarithmic change in the distance to the singularity over that subinterval. Concretely, the number of quadrature nodes in $s$ is the maximum of $q_n/2$ and the greatest integer not exceeding $q_n\log_{10}(h_{max}/h_{min}) / + 5q_n (s_l - s_r)/(2\pi)$. The integration in $s$ is done via Gauss-Legendre quadrature after an exponential transformation. 

\begin{figure}
\includegraphics[width = 0.5\linewidth]{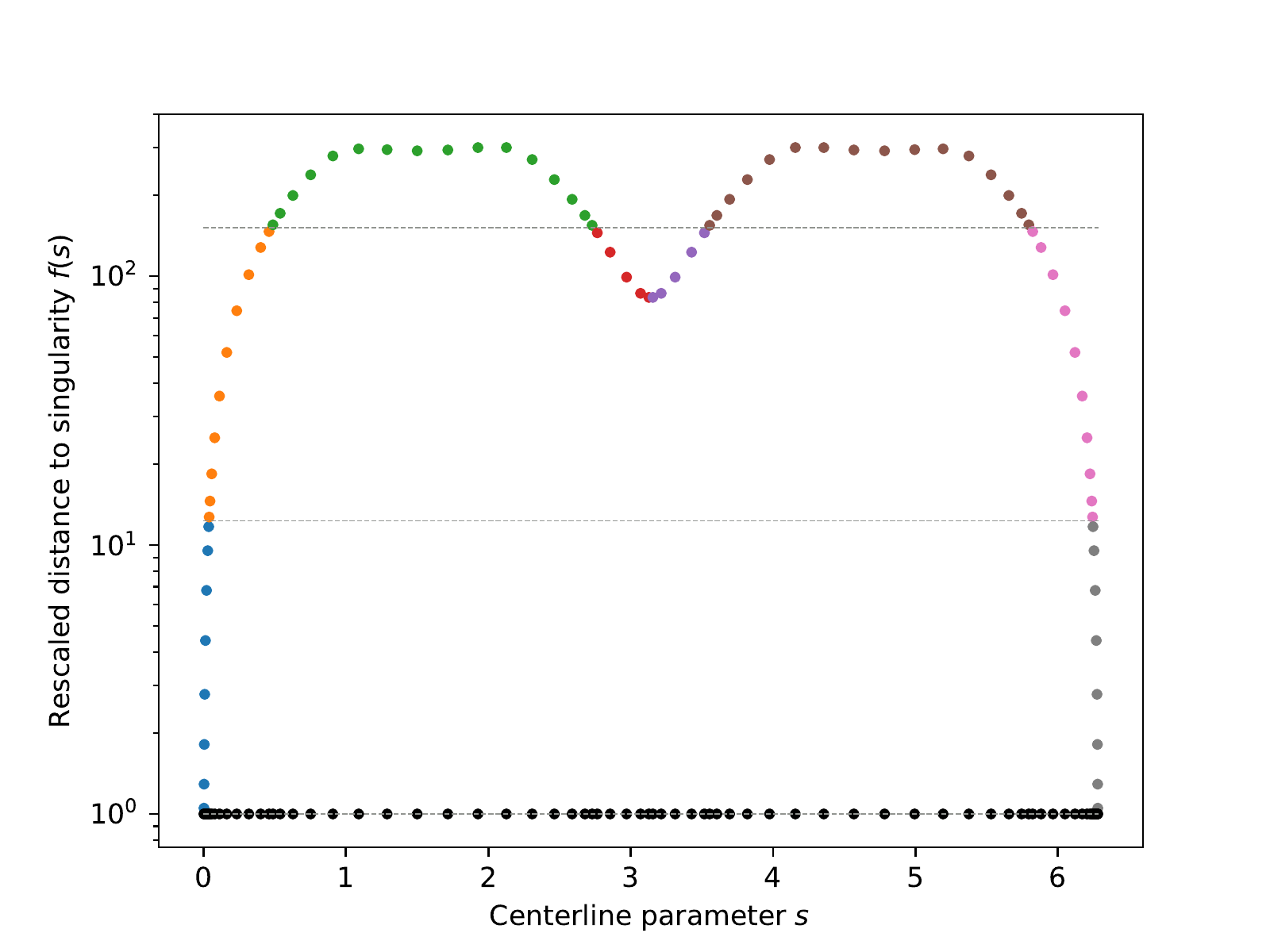}
\caption{If the centerline is a figure-eight shape with a near self-intersection, $\bm\gamma(0) \approx \bm\gamma(\pi)$, we need a quadrature rule that integrates the region near the singularity at $s=0$ and the near-approach at $s=\pi$ carefully. We construct one by splitting the arclength $s$-domain into subregions with boundaries located where the distance to the singularity takes geometrically spaced values, or at extrema where the second derivative is especially large (for example, there are three extrema within the subinterval with green dots, but they are ignored; in contrast, the minimum at $s=\pi$ becomes a subinterval boundary). The number of quadrature points on each subinterval is equal to the quadrature parameter ($q_n=8$ here), multiplied by a factor which increases the number of nodes if the subinterval is large in $s$ or if $f(s)$ takes on very different values inside the interval. 
The logic is complicated but the algorithm runs quickly and allows the use of smaller quadrature rules than would be necessary without this adaptive procedure. }
\label{fig:quadsplits}
\end{figure}
The exponential transformation is as follows. Let the integration interval be $[s_l, s_r]\times [-\pi,\pi]$.  We already know the values $h(s_l)$ and $h(s_r)$ from the subdivision procedure referenced earlier. If these are equal, we use ordinary Gauss-Legendre integration in $s$. If the endpoint values of $h$ are unequal, that is, if one end of the corresponding centerline section is closer to the singularity than the other, we carry out the transformation 
\begin{equation}
\int_{s_l}^{s_r} \psi(s)\,ds = \int_0^1 \psi\Big(A + B\exp(C t)\Big)\cdot BC\exp(Ct)\,dt
\label{eq:sintegralrescaling}
\end{equation}
where 
\begin{equation}
A = \frac{s_r h(s_l) - s_l h(s_r)}{h(s_r)-h(s_l)},\qquad 
B = (s_l-s_r)\frac{h(s_l)}{h(s_r)-h(s_l)},
\qquad 
C = \log(h(s_r) / h(s_l)).
\end{equation}

This transformation comes from the fact that $s(t) =  A + Be^{Ct}$ is the solution of the first-order boundary value problem 
\begin{equation}s(0) = s_\ell,\quad s(1) = s_r,\quad \frac{ds}{dt} = C\left(\frac{s-s_l}{s_r-s_l} h(s_r) + \frac{s_r-s}{s_r-s_l} h(s_l)  \right)\end{equation}
where the expression in large parentheses is the linear interpolation of $h(s)$ on $[s_l,s_r]$, and $C$ is a proportionality constant whose value is determined as part of the solution of the BVP. The consequence is that the product of $ds/dt$ and $1/h(s)$ should be approximately constant, so the transformed problem can be integrated with fewer quadrature nodes. 
We apply Gauss-Legendre integration on the right-hand side of \eqref{eq:sintegralrescaling}. 

\subsection{Convergence of the quadrature procedure}
As the quadrature is an important ingredient in our overall numerical method, and also an independent problem that may have value beyond the present paper, we present some convergence results for three model integrands. These problems have similar behavior to those actually needed for our matrix assembly but we have simplified them somewhat to make the results in this section more replicable by others. 

The integration domains are the tubes whose centerlines are given by \eqref{eq:hairtiecurve} using $H=0.9$, and whose circular cross sections have radius $\epsilon\in\{0.05,0.005,0.0005,0.00005\}$. We do \emph{not} reparameterize or rescale the centerline. 
On each of these domains we consider three integrands which are singular at the point  
\begin{equation}
\bm X(s^*, \theta^*) =  \begin{pmatrix}1/20\\ \sqrt{3}/20 \\0 \end{pmatrix} + \epsilon 
\begin{pmatrix}
1/2\\\sqrt{3}/{2}\\0
\end{pmatrix} .
\end{equation}
This corresponds to taking $s^* = \pi/3$ and $\theta^*=0$ in a system where the centerline frame has $\bm N(\pi/3) = (0.5,0.5\sqrt{3},0)$, that is, we can use the Frenet frame for these tests rather than the more specialized Bishop frame constructed in the previous subsection in order to make the results more accessible.  
Writing $\bm r  = \bm r(s,\theta) = \bm X(s,\theta) - \bm X(s^*,\theta^*)$, the definitions of our three model integrals are: 
\begin{align}
I_{1/r} &= \int_0^{2\pi}\int_0^{2\pi} \frac{1}{r}J(s,\theta)\,d\theta\,ds \\ 
I_{SL} &= \frac{1}{8\pi}\int_0^{2\pi}\int_0^{2\pi} \frac{(\bm r\cdot\hat{\bm x})(\bm r\cdot\hat{\bm y})}{r^3}\cos(50 s + 2\theta)J(s,\theta)\,d\theta\,ds \\ 
I_{DL} &= \frac{-3}{4\pi}\int_0^{2\pi}\int_0^{2\pi} \frac{(\bm r\cdot\hat{\bm x})(\bm r\cdot\hat{\bm y})(\bm r\cdot{\bm \nu})}{r^5}\cos(50 s + 2\theta)J(s,\theta)\,d\theta\,ds.
\end{align}
The first model integrand is the simplest: it is just the reciprocal of the distance to the source point. 
This first integrand has a $1/r$ singularity at the source point, but none of the angular dependence on the divergence rate that the single- and double-layer integrands display. The second integrand is one component of the single-layer potential multiplied by a Fourier basis function.  
The third is one component of the double-layer potential multiplied by a Fourier basis function. 
We have omitted the point-source terms for simplicity and because they are regular. 

As the results of Table \ref{tbl:quadconverge} indicate, we generally get thirteen to fifteen digits for the first integrand using between 10000 and 20000 quadrature nodes. The single- and double-layer integrals $I_{SL}$ and $I_{DL}$ approach zero more quickly as $\epsilon\to0$, but if we count the leading zeros as correct digits then we can claim a similarly good performance. To improve on the last block of the table, that is, to get more than five significant figures in $I_{DL}$ for the smallest value of $\epsilon$, presents an interesting numerical challenge but is not a likely source of significant error in the overall numerical method presented here.  

The quadrature procedure presented here produces good results in the tests we considered, but we hope in the future to find a simpler and more lightweight method.

\begin{table}
\begin{tabular}{cccccc}
$\epsilon$ & $q_n$ & total nodes  & $1/r$ & SL & DL \\\toprule
$5\cdot 10^{-2}$ & $9$ & $2163$ & $7.589782829235038e+00${} & $8.798999514900755e-04${} & $-2.218982313369795e-02${} \\ 
 & $13$ & $4567$ & $7.589782838734336e+00${} & $8.793815847738606e-04${} & $-2.220170156415448e-02${} \\ 
 & $18$ & $9018$ & $7.589782838781613e+00${} & $8.793927729055879e-04${} & $-2.220137549889306e-02${} \\ 
 & $24$ & $15984$ & $7.589782838781604e+00${} & $8.793927729104493e-04${} & $-2.220137549905115e-02${} \\ 
 & $31$ & $26328$ & $7.589782838781603e+00${} & $8.793927729104521e-04${} & $-2.220137549905804e-02${} \\ 
 & $39$ & $41896$ & $7.589782838781612e+00${} & $8.793927729104533e-04${} & $-2.220137549904386e-02${} \\ 
 & $48$ & $63936$ & $7.589782838781607e+00${} & $8.793927729104308e-04${} & $-2.220137549899438e-02${} \\ 
 & $58$ & $93264$ & $7.589782838781607e+00${} & $8.793927729104445e-04${} & $-2.220137549903136e-02${} \\ 
\midrule
$5\cdot 10^{-3}$ & $9$ & $1547$ & $9.256886834463303e-01${} & $6.072425129519274e-06${} & $5.604153771271108e-05${} \\ 
 & $13$ & $3274$ & $9.256898147158956e-01${} & $6.081483386643081e-06${} & $6.352751747862048e-05${} \\ 
 & $18$ & $6390$ & $9.256898207485448e-01${} & $6.082225076874217e-06${} & $6.279410406542996e-05${} \\ 
 & $24$ & $11424$ & $9.256898207525426e-01${} & $6.082225851130320e-06${} & $6.279192290716310e-05${} \\ 
 & $31$ & $18996$ & $9.256898207525437e-01${} & $6.082225851541232e-06${} & $6.279192348132032e-05${} \\ 
 & $39$ & $30138$ & $9.256898207525388e-01${} & $6.082225851527452e-06${} & $6.279192298478882e-05${} \\ 
 & $48$ & $45888$ & $9.256898207525420e-01${} & $6.082225851528601e-06${} & $6.279192339740517e-05${} \\ 
 & $58$ & $67019$ & $9.256898207525425e-01${} & $6.082225851524810e-06${} & $6.279192336884141e-05${} \\ 
\midrule
$5\cdot 10^{-4}$ & $9$ & $1521$ & $1.072470570163802e-01${} & $2.711766253259868e-08${} & $5.397620419211990e-05${} \\ 
 & $13$ & $3211$ & $1.072470604686221e-01${} & $2.712583826063939e-08${} & $5.406876081505809e-05${} \\ 
 & $18$ & $6210$ & $1.072470604378008e-01${} & $2.712525854659307e-08${} & $5.407833085385868e-05${} \\ 
 & $24$ & $11064$ & $1.072470604378060e-01${} & $2.712525879230430e-08${} & $5.407832974464318e-05${} \\ 
 & $31$ & $18538$ & $1.072470604378063e-01${} & $2.712525879697797e-08${} & $5.407833049027967e-05${} \\ 
 & $39$ & $29328$ & $1.072470604378077e-01${} & $2.712525879947280e-08${} & $5.407833147346922e-05${} \\ 
 & $48$ & $44544$ & $1.072470604378010e-01${} & $2.712525876683408e-08${} & $5.407833241615479e-05${} \\ 
 & $58$ & $65018$ & $1.072470604378035e-01${} & $2.712525877772501e-08${} & $5.407833337423479e-05${} \\ 
\midrule
$5\cdot 10^{-5}$ & $9$ & $1719$ & $1.217324007891005e-02${} & $1.828438569716712e-09${} & $5.788310681316210e-06${} \\ 
 & $13$ & $3588$ & $1.217323962076407e-02${} & $1.829372530288228e-09${} & $5.883417669631449e-06${} \\ 
 & $18$ & $6948$ & $1.217323961763665e-02${} & $1.829372728188357e-09${} & $5.883541425090347e-06${} \\ 
 & $24$ & $12432$ & $1.217323961763635e-02${} & $1.829372722135347e-09${} & $5.883533250724377e-06${} \\ 
 & $31$ & $20739$ & $1.217323961763636e-02${} & $1.829372718474100e-09${} & $5.883539417174405e-06${} \\ 
 & $39$ & $32799$ & $1.217323961763877e-02${} & $1.829372719598264e-09${} & $5.883535786441918e-06${} \\ 
 & $48$ & $49776$ & $1.217323961764001e-02${} & $1.829372727088588e-09${} & $5.883559961474092e-06${} \\ 
 & $58$ & $72732$ & $1.217323961763658e-02${} & $1.829372712986068e-09${} & $5.883500032425172e-06${} \\ 
\bottomrule
\end{tabular}
\caption{The numerical integration is challenging because of the singularity and the large aspect ratio of the fiber surface. Here we report convergence results for three model integrands over four orders of magnitude in the fiber radius. We get good results for the simplest problem $I_{1/r}$ using fewer than 10000 function evaluations; for the modified velocity and traction integrals $I_{SL}$ and $I_{DL}$ the accuracy is also good for large $\epsilon$ but decays when $\epsilon\to0$. The results for $I_{DL}$ with $\epsilon=5\cdot10^{-5}$ could particularly be improved, although it is unlikely that this is a significant source of error in our overall algorithm.  }
\label{tbl:quadconverge}
\end{table}

\end{document}